\documentclass{article}
\usepackage{epsfig}
\usepackage{amssymb}
\usepackage{amsmath}
\usepackage{color}
\usepackage{subfigure}
\newtheorem{rem}{Remark}[section]
\newcommand{\br}{\begin{rem}}
\newcommand{\er}{\end{rem}}
\newtheorem{ex}[rem]{Example}
\newcommand{\bex}{\begin{ex}}
\newcommand{\eex}{\end{ex}}
\newtheorem{Def}[rem]{Definition}
\newcommand{\bd}{\begin{Def}}
\newcommand{\ed}{\end{Def}}
\newtheorem{theorem}[rem]{Theorem}
\newcommand{\bt}{\begin{theorem}}
\newcommand{\et}{\end{theorem}}
\newtheorem{prop}[rem]{Proposition}
\newcommand{\bp}{\begin{prop}}
\newcommand{\ep}{\end{prop}}
\newtheorem{lemma}[rem]{Lemma}
\newcommand{\bl}{\begin{lemma}}
\newcommand{\el}{\end{lemma}}
\newtheorem{corollary}[rem]{Corollary}
\newcommand{\bc}{\begin{corollary}}
\newcommand{\ec}{\end{corollary}}
\newcommand{\be}{\begin{equation}}
\newcommand{\ee}{\end{equation}}
\newcommand{\bea}{\begin{eqnarray}}
\newcommand{\eea}{\end{eqnarray}}

\newcommand{\nn}{\nonumber}
\newcommand{\adots}{\mathinner{\mkern2mu\raise1pt\hbox{.}\mkern2mu
\raise4pt\hbox{.}\mkern2mu\raise7pt\hbox{.}\mkern1mu}}

\title{Cluster Mutation-Periodic Quivers and Associated Laurent Sequences}
\author{Allan P. Fordy and Robert J. Marsh,
  \\ School of Mathematics, \\
University of Leeds. \\ Leeds LS2 9JT, UK.\\
e-mail: allan@maths.leeds.ac.uk and marsh@maths.leeds.ac.uk}

\begin{document}

\date{Revised version August 2010}
\maketitle
\bibliographystyle{plain}

\begin{abstract}
We consider quivers/skew-symmetric matrices under the action of mutation (in
the cluster algebra sense). We classify those which are isomorphic to their own
mutation via a cycle permuting all the vertices, and give families of quivers
which have higher periodicity.

The periodicity means that sequences given by recurrence relations arise in a
natural way from the associated cluster algebras.  We present a number of
interesting new families of nonlinear recurrences, necessarily with the Laurent
property, of both the real line and the plane, containing integrable maps as
special cases.  In particular, we show that some of these recurrences can be
linearised and, with certain initial conditions, give integer sequences which
contain all solutions of some particular Pell equations.  We extend our
construction to include recurrences with {\em parameters}, giving an
explanation of some observations made by Gale.

Finally, we point out a connection between quivers which arise in our
classification and those arising in the context of quiver gauge theories.
\end{abstract}

\section{Introduction}

\begin{figure}[htb]
\centering
\subfigure[The Somos $4$ quiver, $S_4$.]{
\includegraphics[width=3.5cm]{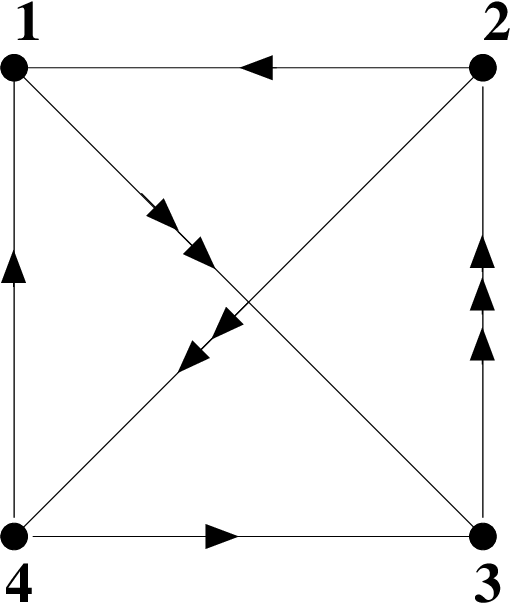}\label{subfig:somos4quiver}
} \qquad\qquad\qquad
\subfigure[Mutation of $S_4$ at $1$.]{
\includegraphics[width=3.5cm]{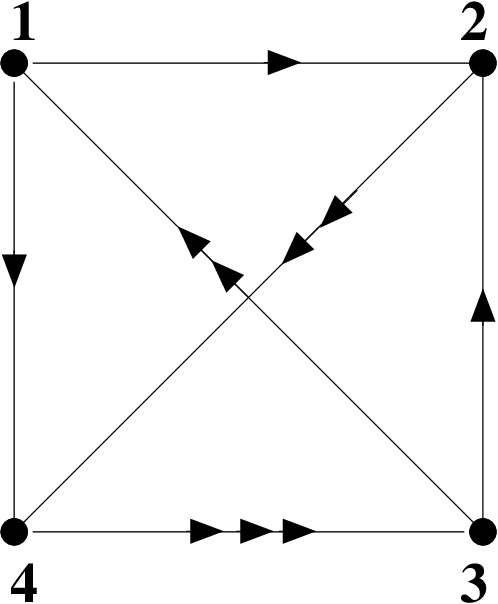}\label{subfig:somos4quiverb}
} \caption{The Somos $4$ quiver and its mutation at $1$.}\label{s4quiver}
\end{figure}

Our main motivation for this work is the connection between cluster algebras
and integer sequences which are Laurent polynomials in their initial
terms~\cite{02-2}. A key example of this is the Somos 4 sequence, which is
given by the following recurrence:
\be  \label{somos4}  %
x_n x_{n+4} = x_{n+1}x_{n+3}+x_{n+2}^2.
\ee  %
This formula, with appropriate relabelling of the variables, coincides with the
cluster exchange relation~\cite{02-3} (recalled below; see
Section~\ref{laurent}) associated with the vertex $1$ in the
quiver $S_4$ of Figure \ref{subfig:somos4quiver}.
Mutation of $S_4$ at $1$ (as in~\cite{02-3}; see Definition~\ref{d:mutate} below)
gives the quiver shown in
Figure~\ref{subfig:somos4quiverb} and transforms the cluster
$(x_1,x_2,x_3,x_4)$ into $(\tilde x_1,x_2,x_3,x_4)$, where $\tilde x_1$ is
given by
$$
x_1 \tilde x_1 = x_2 x_4+x_3^2 .
$$
Remarkably, after this {\em complicated operation} of mutation on the quiver,
the result is a simple rotation, corresponding to the relabelling of indices
$(1,2,3,4)\mapsto (4,1,2,3)$.  Therefore, a mutation of the new quiver at $2$
gives the {\em same formula} for the exchange relation (up to a relabelling).
It is this simple property that allows us to think of an infinite sequence of
such mutations as iteration of recurrence~(\ref{somos4}).

In this paper, we classify quivers with this property. In this way we obtain a
classification of maps which could be said to be `of Somos type'. In fact we
consider a more general type of ``mutation periodicity'', which corresponds to
Somos type sequences of higher dimensional spaces.

It is interesting to note that many of the quivers which have occurred in the
theoretical physics literature concerning {\em supersymmetric quiver gauge
theories} are particular examples from our classification; see for
example~\cite[\S4]{01-10}. We speculate that some of our other examples may be
of interest in that context.

We now describe the contents of the article in more detail. In Section
\ref{per-prop}, we recall matrix and quiver mutation from \cite{02-3}, and
introduce the notion of periodicity we are considering. It turns out to be
easier to classify periodic quivers if we assume that certain vertices are
sinks; we call such quivers \emph{sink-type}. In Section \ref{p1}, we classify
the sink-type quivers of period 1 as nonnegative integer combinations of a
family of \emph{primitive} quivers. In Section \ref{p2}, we do the same for
sink-type period $2$ quivers, and in Section \ref{p-high} we classify the
sink-type quivers of arbitrary period.

In Section \ref{p1-gen}, we give a complete classification of all period $1$
quivers (without the sink assumption), and give some examples. It turns out the
arbitrary period $1$ quivers can be described in terms of the primitives with
$N$ nodes, together with the primitives for quivers with $N'$ nodes for all $N'$
less than $N$ of the same parity (Theorem \ref{p1-theorem}).

In Section \ref{p2-quivers}, we classify quivers of period $2$ with at most $5$
nodes. These descriptions indicate that a full classification for higher period
is likely to be significantly more complex than the classification of period
$1$ quivers.  However, it is possible to construct a large family of period $2$
(not of sink-type) quivers, which we present in Section \ref{p2-reg}.

In Section \ref{laurent}, we describe the recurrences that can be associated to
period 1 and period 2 quivers via Fomin-Zelevinsky cluster mutation. The nature
of the cluster exchange relation means that the recurrences we have associated
to periodic quivers are in general nonlinear. However, in Section \ref{linear},
we show that the recurrences associated to period 1 primitives can be
linearised. This allows us to conclude in Section \ref{linear} that certain
simple linear combinations of subsequences of the first primitive period 1
quiver (for arbitrarily many nodes) provide all the solutions to an associated
Pell equation.

In Section \ref{s:coefficients} we extend our construction of mutation periodic
quivers to include quivers with {\em frozen} cluster variables, thus enabling
the introduction of {\em parameters} into the corresponding recurrences.  As a
result, we give an explanation of some observations made by Gale in
\cite{91-17}.

In Section \ref{susy}, we give an indication of the connections with
supersymmetric quiver gauge theories.

In Section \ref{conclude}, we present our final conclusions. Section
\ref{append} is an appendix to Section \ref{linear}.

\section{The Periodicity Property}\label{per-prop}

We consider quivers with no $1$ cycles or $2$-cycles (i.e.\ the quivers
on which cluster mutation is defined). Any $1$- or $2$-cycles which
arise through operations on the quiver will be cancelled.
The vertices of $Q$ will be assumed to lie on the vertices of a regular
$N$-sided polygon, labelled $1,2,\ldots ,N$ in clockwise order.

In the usual way, we shall identify a quiver $Q$, with $N$ nodes, with the
unique skew-symmetric $N\times N$ matrix $B_Q$ with $(B_{Q})_{ij}$ given by the
number of arrows from $i$ to $j$ minus the number of arrows from $j$ to $i$..
We next recall the definition of quiver mutation~\cite{02-3}.

\bd[Quiver Mutation]\label{d:mutate}
Given a quiver $Q$ we can mutate at any of its nodes.
The mutation of $Q$ at node $k$, denoted by $\mu_k Q$, is constructed (from
$Q$) as follows:
\begin{enumerate}
\item  Reverse all arrows which either originate or terminate at node $k$.
\item  Suppose that there are $p$ arrows from node $i$ to node $k$ and $q$
arrows from node $k$ to node $j$ (in $Q$).
Add $pq$ arrows going from node $i$ to node $j$ to any arrows already there.
\item Remove (both arrows of) any two-cycles created in the previous
steps.
\end{enumerate}
\ed    %
Note that Step $3$ is independent of any choices made in the removal
of the two-cycles, since the arrows are not labelled.
We also note that in Step $2$, $pq$ is just the number of paths of length $2$
between nodes $i$ and $j$ which pass through node $k$.

\br[Matrix Mutation]   %
Let $B$ and $\tilde B$ be the skew-symmetric matrices corresponding to the
quivers $Q$ and $\tilde Q=\mu_k Q$.  Let $b_{ij}$ and $\tilde b_{ij}$ be the
corresponding matrix entries.  Then quiver mutation amounts to the following
formula
\be   \label{gen-mut}  %
\tilde b_{ij}= \left\{ \begin{array}{ll}
                        -b_{ij} & \mbox{if}\;\; i=k\;\;\mbox{or}\;\; j=k, \\
                        b_{ij}+\frac{1}{2} (|b_{ik}|b_{kj}+b_{ik}|b_{kj}|)
                        & otherwise.
                        \end{array}  \right.
\ee    %
\er      %

This is the original formula appearing (in a more general context)
in~\cite{02-3}.

We number the nodes from $1$ to $N$, arranging them equally spaced on a circle
(clockwise ascending). We consider the permutation
$\rho:(1,2,\cdots ,N)\rightarrow (N,1,\cdots ,N-1)$.
Such a permutation acts on a quiver $Q$ in such a way that the number of arrows
from $i$ to $j$ in $Q$ is the same as the number of arrows from $\rho^{-1}(i)$
to $\rho^{-1}(j)$ in $\rho Q$. Thus the arrows of $Q$ are rotated clockwise
while the nodes remain fixed (alternatively, this operation can be interpreted
as leaving the arrows fixed whilst the nodes are moved in an anticlockwise
direction). We will always fix the positions of the nodes in our diagrams.

Note that the action $Q\mapsto \rho Q$ corresponds to the conjugation
$B_Q\mapsto \rho B_Q \rho^{-1}$, where
$$
\rho =  \left( \begin{array}{cccc}
                   0 & \cdots & \cdots & 1 \\
                  1 & 0 & & \vdots \\
                     &\ddots & \ddots & \vdots \\
                      & & 1 & 0
                      \end{array} \right)
$$
(we will use the notation $\rho$ for both the permutation and corresponding matrix).

We consider a sequence of mutations, starting at node $1$,
followed by node $2$, and so on. Mutation at node $1$ of a quiver $Q(1)$ will
produce a second quiver $Q(2)$.  The mutation at node $2$ will therefore be of
quiver $Q(2)$, giving rise to quiver $Q(3)$ and so on.

\bd
We will say that a quiver $Q$ has \textit{period $m$} if it satisfies
$Q(m+1)=\rho^m Q(1)$, with the mutation sequence depicted by
\be   \label{periodchain}   %
Q=Q(1) \stackrel{\mu_1}{\longrightarrow} Q(2) \stackrel{\mu_2}{\longrightarrow}
  \cdots \stackrel{\mu_{m-1}}{\longrightarrow} Q(m)
       \stackrel{\mu_m}{\longrightarrow} Q(m+1)=\rho^m Q(1).
\ee    %
We call the the above sequence of quivers the \emph{periodic chain}
associated to $Q$.
\ed

Note that permutations other than $\rho^m$ could be used here, but we
do not consider them in this article. If $m$ is minimal in the above, we say
that $Q$ is strictly of period $m$. Also note that each of the quivers
$Q(1),\ldots,Q(m)$ is of period $m$ (with a renumbering of the vertices),
if $Q$ is.

Recall that a node $i$ of a quiver $Q$ is said to be a \emph{sink}
if all arrows incident with $i$ end at $i$, and is said to be a
\emph{source} if all arrows incident with $i$ start at $i$.

\br[Admissible sequences] \label{r:admissible}
Recall that an \emph{admissible sequence} of sinks in
an acyclic quiver Q is a total ordering $v_1,v_2,\ldots ,v_N$ of its vertices
such that $v_1$ is a sink in $Q$ and $v_i$ is a sink in
$\mu_{v_{i-1}}\mu_{v_{i-2}}\cdots \mu_{v_1}(Q)$ for $i=2,3,\ldots ,N$. Such a
sequence always has the property that $\mu_{v_N}\mu_{v_{N-1}}\cdots
\mu_{v_1}(Q)=Q$ \cite[\S5.1]{06-2}. This notion is of importance in the
representation theory of quivers.

We note that if any (not necessarily acyclic) quiver $Q$ has period 1 in our
sense, then $\mu_1 Q=\rho Q$. It follows that $\mu_N \mu_{N-1} \cdots \mu_1
Q=Q$. Thus any period $1$ quiver has a property which can be regarded as a
generalisation of the notion of existence of an admissible sequence of sinks.
In fact, higher period quivers also possess this property provided the period
divides the number of vertices.
\er  %

\section{Period $1$ Quivers}\label{p1}

We now introduce a finite set of particularly simple quivers of period $1$,
which we shall call the \emph{period $1$ primitives}.
Remarkably, it will later be seen that in a certain sense they form a
``basis'' for the set of all quivers of period $1$.
We shall also later see that period $m$ primitives can be defined as certain
sub-quivers of the period $1$ primitives.

\bd[Period $1$ sink-type quivers]\label{d:period1}
A quiver $Q$ is said to be a \emph{period $1$ sink-type quiver}
if it is of period $1$ and node $1$ of $Q$ is a sink.
\ed

\bd[Skew-rotation] \label{d:skewrotation}
We shall refer to the matrix
$$\tau = \left(\begin{array}{cccc}
                   0 & \cdots & \cdots & -1 \\
                  1 & 0 & & \vdots \\
                     &\ddots & \ddots & \vdots \\
                      & & 1 & 0
                      \end{array} \right).
$$
as a {\em skew-rotation}.
\ed

\bl[Period $1$ sink-type equation]\label{l:period1sink}
A quiver $Q$ with a sink at $1$ is period $1$ if and only if
$\tau B_Q \tau^{-1}=B_Q$.
\el

\noindent \textbf{Proof}:
If node $1$ of $Q$ is a sink,
there are no paths of length $2$ through it, and the second part of Definition
\ref{d:mutate} is void.  Reversal of the arrows at node $1$ can be done
through a simple conjugation of the matrix $B_Q$:
$$
\mu_1 B_Q = D_1 B_Q D_1, \quad\mbox{where}\quad D_1=
\mbox{diag}(-1,1,\cdots ,1).
$$
Equating this to $\rho B_Q \rho^{-1}$ leads to the equation
$\tau B_Q \tau^{-1} = B_Q$ as required, noting that
$$\tau=D_1\rho.$$
$\Box$

The map $M\mapsto \tau M \tau^{-1}$ simultaneously cyclically permutes the
rows and columns of $M$ (up to a sign), while $\tau^N=-I_N$,
hence $\tau$ has order $N$.
This gives us a method for building period $1$ matrices: we sum over
$\tau$-orbits.

\paragraph{The period $1$ primitives $P_N^{(k)}$.}
We consider a quiver with just a single arrow from $N-k+1$ to $1$,
represented by the skew-symmetric matrix $R_N^{(k)}$ with
$(R_N^{(k)})_{N-k+1,1}=1, (R_N^{(k)})_{1,N-k+1}=-1$ and
$(R_N^{(k)})_{ij}=0$ otherwise.

We define skew-symmetric matrices $B_N^{(k)}$ as follows:
\begin{equation} \label{e:primitives}
B_N^{(k)} = \left\{ \begin{array}{ll}
  \sum_{i=0}^{N-1} \tau^i R_N^{(k)} \tau^{-i}, &
       \mbox{if\ }N=2r+1\mbox{\ and\ }1\leq k\leq r,
            \mbox{\ or\ if\ } N=2r\mbox{\ and\ }1\leq k\leq r-1; \\
\sum_{i=0}^{r-1} \tau^i R_N^{(r)} \tau^{-i}, &
                   \mbox{if\ }k=r \mbox{\ and\ }N=2r. \\
                \end{array}\right.
\end{equation}
Let $P_N^{(k)}$ denote the quiver corresponding to $B_N^{(k)}$.
We remark that the geometric action of $\tau$ in the above sum is to rotate
the arrow clockwise without change of orientation, except that when the tail
of the arrow ends up at node $1$ it is reversed. It follows that $1$
is a sink in the resulting quiver. Since it is a sum over a $\tau$-orbit,
we have $\tau B_N^{(k)} \tau^{-1} = B_N^{(k)}$, and thus that $P_N^{(k)}$
is a period $1$ sink-type quiver. In fact, we have the simple description:
$$
B_N^{(k)}= \left\{ \begin{array}{ll}
              \tau^k-(\tau^{t})^k, & \mbox{if\ }N=2r+1
                 \mbox{\ and\ }1\leq k\leq r, \mbox{\ or\ }
                           N=2r\mbox{\ and\ }1\leq k\leq r-1; \\
                \tau^r , & \mbox{if\ }N=2r\mbox{\ and\ } k=r.
                \end{array}\right.
$$
where $\tau^{t}$ denotes the transpose of $\tau$.

Note that we have restricted to the choice $1\leq k\leq r$
because when $k>r$, our construction gives nothing new. Firstly,
consider the case $N\not=2k$. Then
$B_N^{(N+1-k)}=B_N^{(k)}$, because the primitive $B_N^{(k)}$ has
exactly two arrows ending at $1$: those starting at $k+1$ and at
$N+k-1$. Starting with either of these arrows produces the same result.
If $N=2k$, these two arrows are identical, and since
$\tau^k$ is skew-symmetric, $\tau^k-(\tau^t)^k=2\tau^k$.  The sum over
$N=2k$ terms just goes twice over the sum over $k$ terms.

In this construction we could equally well have chosen node $1$ to be a source.
We would then have $R_N^{(k)}\mapsto -R_N^{(k)}$, $B_N^{(k)}\mapsto -B_N^{(k)}$
and $P_N^{(k)}\mapsto (P_N^{(k)})^{opp}$, where $Q^{opp}$ denotes the
opposite quiver of $Q$ (with all arrows reversed).
Our original motivation in terms of sequences with the Laurent property
is derived through cluster exchange relations, which do not distinguish
between a quiver and its opposite, so we consider these as equivalent.

\br  %
We note that each primitive is a disjoint union of cycles or arrows, i.e.
quivers whose underlying graph is a union of components which are either of
type $A_2$ or of type $\widetilde{A}_m$ for some $m$.
\er  %

\begin{figure}[htb]
\centering
\subfigure[$P_2^{(1)}$]{
\includegraphics[width=2.5cm]{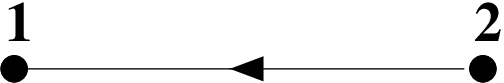}\label{subfig:P21}
}\qquad
\subfigure[$P_3^{(1)}$]{
\includegraphics[width=2.5cm]{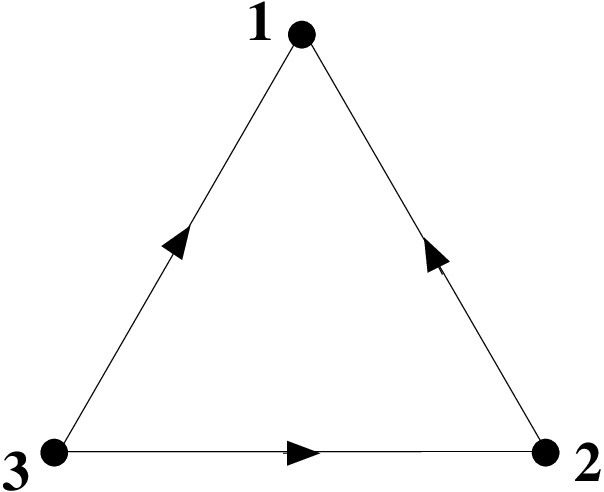}\label{subfig:P31}
}\qquad
\subfigure[$P_4^{(1)}$]{
\includegraphics[width=2.5cm]{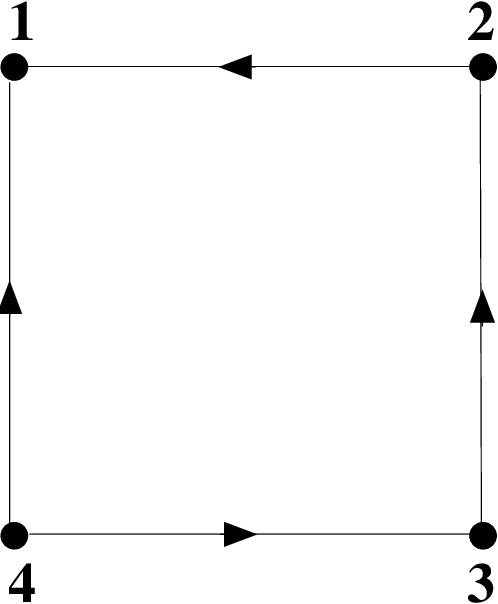}\label{subfig:P41}
}\qquad
\subfigure[$P_4^{(2)}$]{
\includegraphics[width=2.5cm]{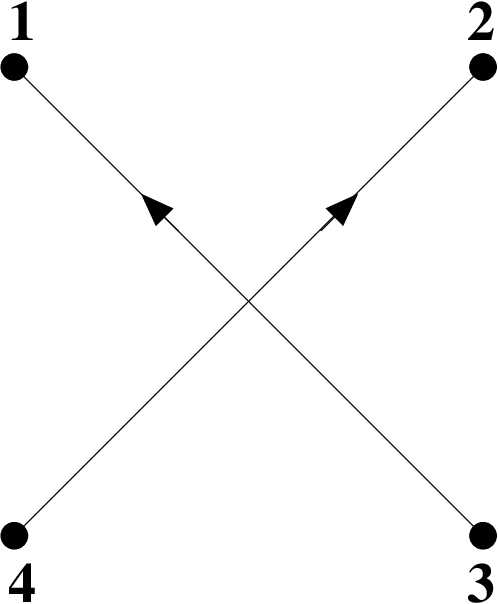}\label{subfig:P42}
}
\caption{The period $1$ primitives for $2$, $3$ and $4$ nodes.}
\label{234node}
\end{figure}

\begin{figure}[htb]
\centering
\subfigure[$P_5^{(1)}$]{
\includegraphics[width=3cm]{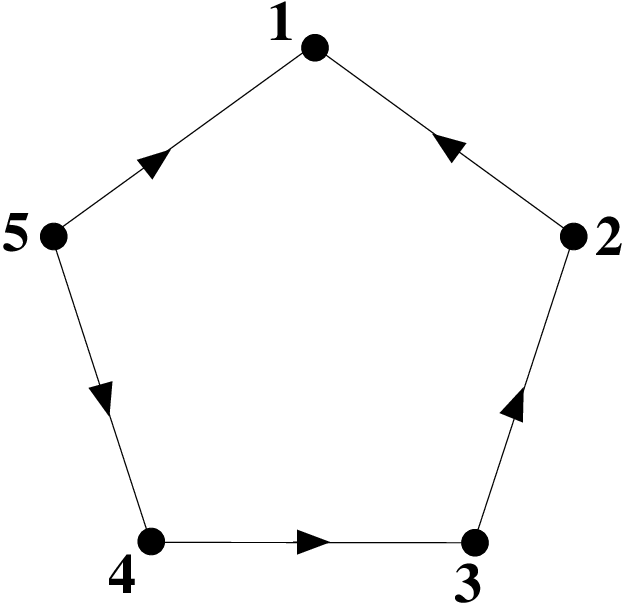}\label{subfig:P51}
}\qquad
\subfigure[$P_5^{(2)}$]{
\includegraphics[width=3cm]{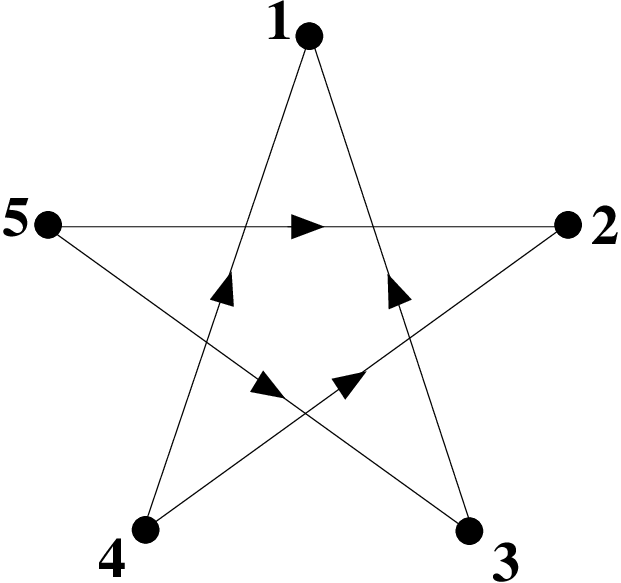}\label{subfig:P52}
}
\caption{The period $1$ primitives for $5$ nodes.}
\label{5node}
\end{figure}

\begin{figure}[htb]
\centering
\subfigure[$P_6^{(1)}$]{
\includegraphics[width=3cm]{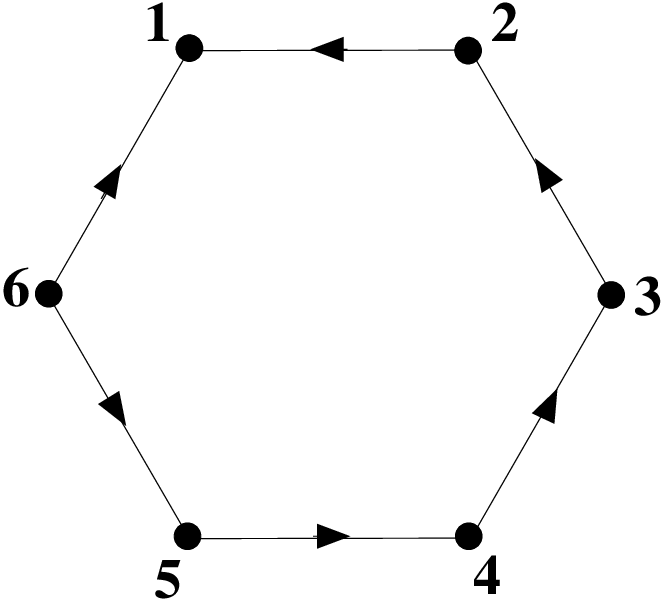}\label{subfig:P61}
}\qquad
\subfigure[$P_6^{(2)}$]{
\includegraphics[width=3cm]{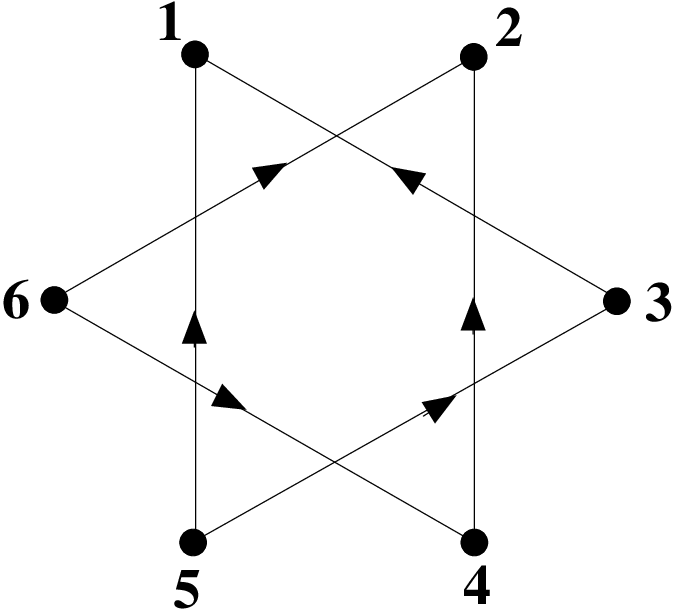}\label{subfig:P62}
}\qquad
\subfigure[$P_6^{(3)}$]{
\includegraphics[width=3cm]{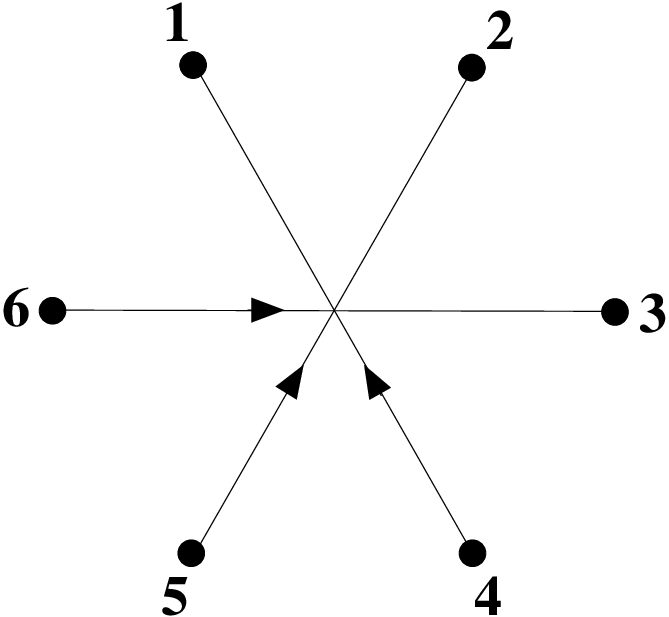}\label{subfig:P63}
}
\caption{The period $1$ primitives for $6$ nodes.}
\label{6node}
\end{figure}

Figures~\ref{234node} to \ref{6node} show the period $1$ primitives we have
constructed, for $2\leq N\leq 6$.

\br[An Involution $\iota : Q\mapsto Q^{opp}$] \label{iota} %
It is easily seen that the following permutation of the nodes is a symmetry of
the primitives $P_N^{(i)}$ (if we consider $Q$ and $Q^{opp}$ as equivalent):
$$
\iota:(1,2,\cdots ,N)\mapsto (N,N-1,\cdots ,1) .
$$
In matrix language, this follows from the facts that $\iota R_N^{(k)}\iota =
-R_N^{(k)}$ and $\iota \tau\iota = \tau^t$, where $\iota =
\left(\begin{array}{cccc}
                   0 &  &  & 1 \\
                   &  & \adots &  \\
                     & \adots &  &  \\
                     1 & &  & 0
                      \end{array} \right).
$

It is interesting to note that $\rho$ is a Coxeter element in $\Sigma_N$
regarded as a Coxeter group, while $\iota$ is the longest element.
\er   %

We may combine primitives to form more complicated quivers.  Consider the
sum
$$
P = \sum_{i=1}^r m_i P_N^{(i)},
$$
where $N=2r$ or $2r+1$ for $r$ an integer and the $m_i$ are arbitrary integers.
It is easy to see that the corresponding quiver is a period $1$ sink-type
quiver whenever $m_i\geq 0$ for all $i$. In fact, we have:

\bp[Classification of period $1$ sink-type quivers]\label{p1sinktypeclassification}
Let $N=2r$ or $2r+1$, where $r$ is an integer.
Every period $1$ sink-type quiver with $N$ nodes has corresponding matrix of
the form  $B = \sum_{k=1}^r m_k B_N^{(k)}$, where the $m_k$ are arbitrary
nonnegative integers.
\ep

\noindent {\bf Proof}: Let $B$ be the matrix of a period $1$ sink-type quiver.
It remains to show that $B$ is of the form stated. We note that conjugation by
$\tau$ permutes the set of summands appearing in the
definition~(\ref{e:primitives}) of the
$B_N^{(k)}$, i.e.\ the elements $\tau^i R_N^{(k)} \tau^{-i}$ for $0\leq i\leq
N-1$ and $1\leq k\leq r$ if $N=2r+1$, for $0\leq i \leq N-1$ and $1\leq k\leq
r-1$ if $N=2r$, together with the elements  $\tau^i R_N^{(r)} \tau^{-i}$ for
$0\leq i\leq r-1$ if $N=2r$. These $\frac{1}{2}N(N-1)$ elements are easily seen
to form a basis of the space of real skew-symmetric matrices. By
Lemma~\ref{l:period1sink}, $\tau B \tau^{-1}=B$, so $B$ is a linear combination
of the period $1$-primitives (which are the orbit sums for the conjugation
action of $\tau$ on the above basis), $B = \sum_{k=1}^r m_k B_N^{(k)}$.
The support of the $B_N^{(k)}$ for $1\leq k\leq r$ is distinct, so
$B_{N-k+1,1}=m_k$ for $1\leq k\leq r$ (where the support of a matrix is
the set of positions of its non-zero entries).
Hence the $m_k$ are integers, as $B$ is an integer matrix. Since $B$ is sink-type,
all the $m_k$ must be nonnegative.$\Box$

Note that this means all period $1$ sink-type quivers are invariant under
$\iota$ in the above sense.
We also note that if the $m_k$ are taken to be of mixed sign, then $Q$ is no
longer periodic without the addition of further ``correction'' terms.
Theorem \ref{p1-theorem} gives these correction terms.

\section{Period $2$ Quivers}\label{p2}

Period $2$ primitives will be defined in a similar way. First, we make
the following definition:

\bd[Period $2$ sink-type quivers]\label{d:period2}
A quiver $Q$ is said to be a \emph{period $2$ sink-type quiver}
if it is of period $2$, node $1$ of $Q(1)=Q$ is a sink, and node $2$ of
$Q(2)=\mu_1 Q$ is a sink.
\ed

Let $Q$ be a period $2$ quiver. Then we have two quivers in our periodic
chain (\ref{periodchain}), $Q(1)$ and $Q(2)=\mu_1 Q$,
with corresponding matrices $B(1), B(2)$.
If $Q(1)$ is of sink-type then, since node $1$ is a sink in $Q(1)$,
the mutation $Q(1)\mapsto \mu_1 Q(1)=Q(2)$ again only involves the reversal of
arrows at node $1$. Similarly, since node $2$ is a sink for $Q(2)$,
the mutation $Q(2)\mapsto \mu_2 Q(2)$ only involves the reversal of arrows
at node $2$.

Obviously each period $1$ quiver $Q$ is also period $2$, where $B(2) = \rho
B(1) \rho^{-1}$. However, we will construct some strictly period $2$
primitives.

As before, we have:

\bl[Period $2$ sink-type equation]
Suppose that $Q$ is a quiver with a sink at $1$ and that $Q(2)$
has a sink at $2$. Then $Q$ is period $2$ if and only if
$\tau^2 B_Q \tau^{-2}=B_Q$.
\el

\noindent {\bf Proof}:
As before, reversal of the arrows at node $1$ of $Q$ can be achieved through
a simple conjugation of its matrix: $\mu_1 B_Q = D_1 B_Q D_1$.
Similarly, reversal of the arrows at node $2$ of $Q(2)$ can be achieved through
$$
\mu_2 B_Q(2) = D_2 B_Q(2) D_2, \quad\mbox{where}\quad D_2=
\mbox{diag}(1,-1,1,\cdots ,1) = \rho D_1 \rho^{-1} .
$$
Equating the composition to $\rho^2 B_Q \rho^{-2}$
leads to the equation
$$
B_Q = D_1D_2 \rho^2 B_Q \rho^{-2}D_2D_1 =
         \tau^2 B_Q \tau^{-2}.
$$
$\Box$

Following the same procedure as for period $1$, we need to form orbit-sums
for $\tau^2$ on the basis considered in the previous section; we shall call
these period $2$ primitives.

A $\tau$-orbit of odd cardinality is also a $\tau^2$-orbit, so the orbit
sum will be a period $2$ primitive which is also of period $1$. Thus we
cannot hope to get period $2$ solutions which are not also period $1$
solutions unless there are an even number of nodes.
A $\tau$-orbit of even cardinality splits into two $\tau^2$-orbits.

When $N=2r$, the matrices $R_N^{(k)}$, for $1\leq k\leq r-1$, generate
strictly period $2$ primitives $P_{N,2}^{(k,1)}$, with matrices given by
$$
B_{N,2}^{(k,1)} = \sum_{i=0}^{r-1} \tau^{2i} R_N^{(k)} \tau^{-2i},
$$
If, in addition, $N$ is divisible by $4$, we obtain the additional strictly
period $2$ primitives $P_{N,2}^{(r,1)}$, with matrices given by:
$$
B_{N,2}^{(r,1)} = \sum_{i=0}^{r/2-1} \tau^{2i} R_N^{(r)} \tau^{-2i}.
$$

Geometrically, the primitive $P_{N,2}^{(k,1)}$ is obtained from the period $1$
primitive $P_N^{(k)}$ by ``removing half the arrows'' (the ones corresponding
to odd powers of $\tau$). The removed arrows form another period $2$ primitive,
called $P_{N,2}^{(k,2)}$, which may be defined as the matrix:
$$
B_{N,2}^{(k,2)} =  \tau B_{N,2}^{(k,1)} \tau^{-1}.
$$

We make the following observation:

\bl \label{relation} %
We have $$\rho^{-1}\mu_1 B_{N,2}^{(k,1)} \rho = B_{N,2}^{(k,2)}$$
for $1\leq k\leq r$.
\el   %

\noindent \textbf{Proof}:
For $1\leq k\leq r-1$, we have
\begin{eqnarray*}
\rho^{-1} \mu_1 B_{N,2}^{(k,1)} \rho & = &
\rho^{-1} D_1 B_{N,2}^{(k,1)} D_1^{-1} \rho \\
&=& \tau^{-1} \left(\sum_{i=0}^{r-1} \tau^{2i} R_N^{(k)} \tau^{-2i}\right)\tau \\
&=& \tau \left(\sum_{i=0}^{r-1} \tau^{2i-2} R_N^{(k)}
\tau^{2-2i}\right)\tau^{-1},
\end{eqnarray*}
since $\rho^{-1}D_1=\tau^{-1}$.  Since $\tau^{-2}=-\tau^{2r-2}$, we have
$\rho^{-1} B_{N,2}^{(k,1)}(2) \rho  = \tau B_{N,1}^{(k,1)} \tau^{-1} =
B_{N,2}^{(k,2)}$.
A similar argument holds for $k=r$, noting that in this case,
$\tau^{-2}R_N^{(k)}=\tau^{r-2}R_N^{(k)}$.
$\Box$

\begin{figure}[htb]
\centering
\subfigure[$P_{4,2}^{(1,1)}$]{
\includegraphics[width=2.5cm]{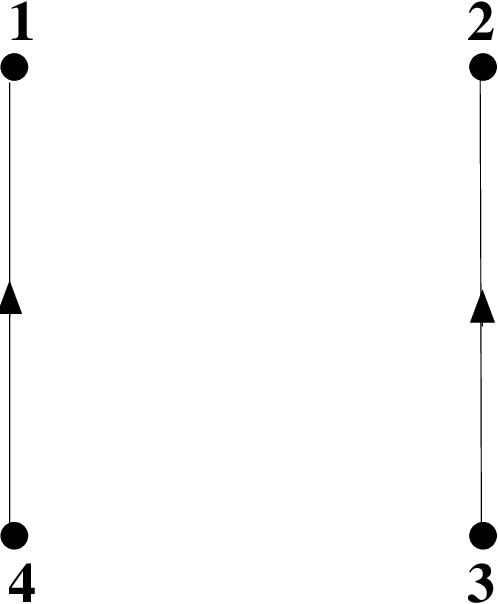}\label{subfig:P4211}
}\qquad
\subfigure[$P_{4,2}^{(1,2)}$]{
\includegraphics[width=2.5cm]{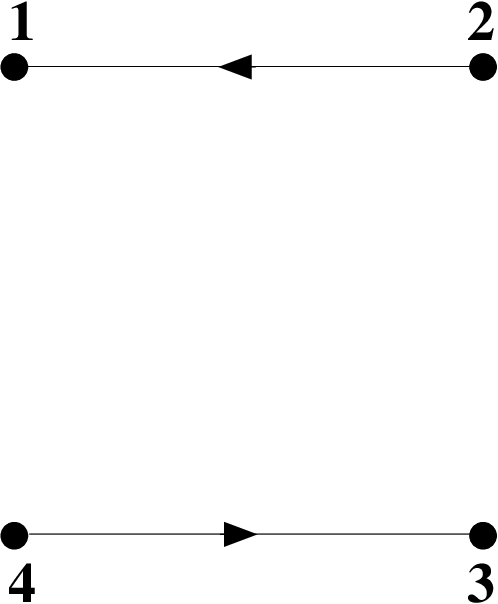}\label{subfig:P4212}
}\qquad
\subfigure[$P_{4,2}^{(2,1)}$]{
\includegraphics[width=2.5cm]{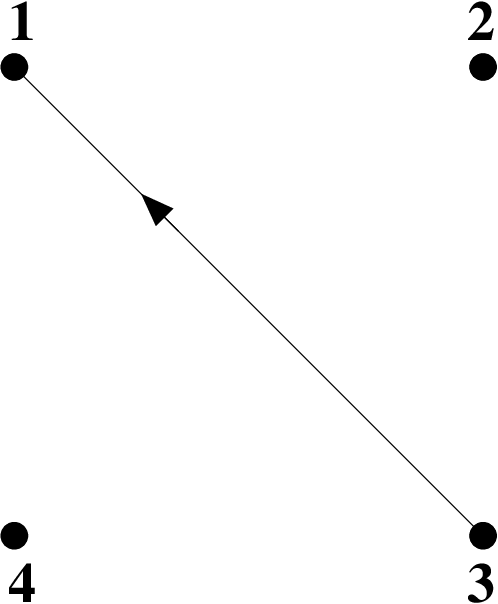}\label{subfig:P4221}
}\qquad
\subfigure[$P_{4,2}^{(2,2)}$]{
\includegraphics[width=2.5cm]{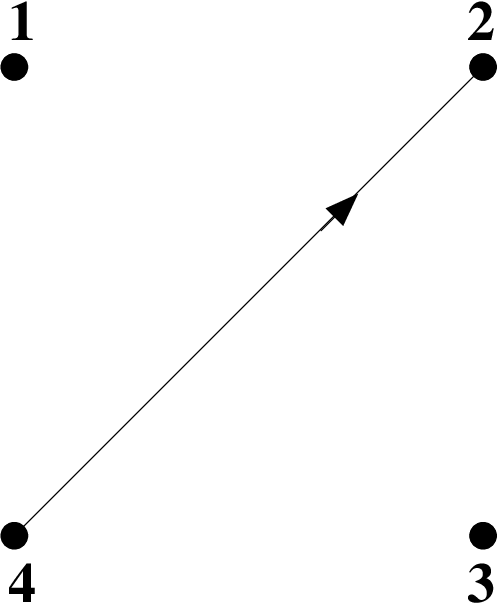}\label{subfig:P4222}
}
\caption{The strictly period $2$ primitives for $4$ nodes.}
\label{4node-p2}
\end{figure}

\begin{figure}[htb]
\centering
\subfigure[$P_{6,2}^{(1,1)}$]{
\includegraphics[width=3cm]{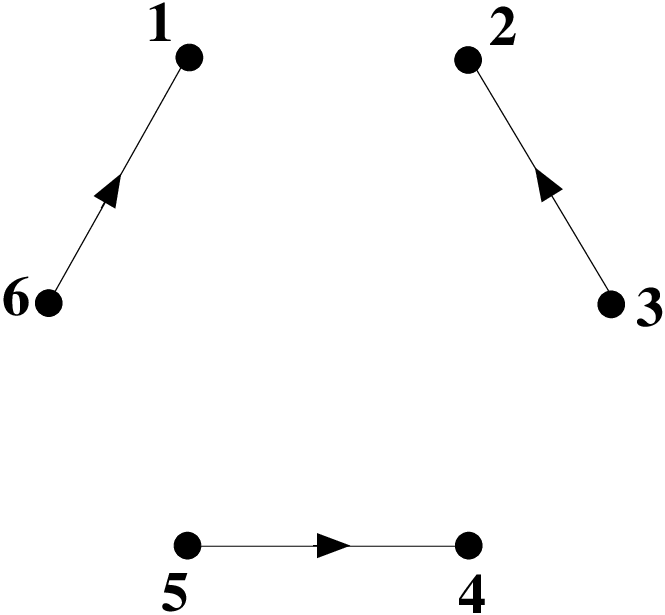}\label{subfig:P6211}
}\qquad
\subfigure[$P_{6,2}^{(1,2)}$]{
\includegraphics[width=3cm]{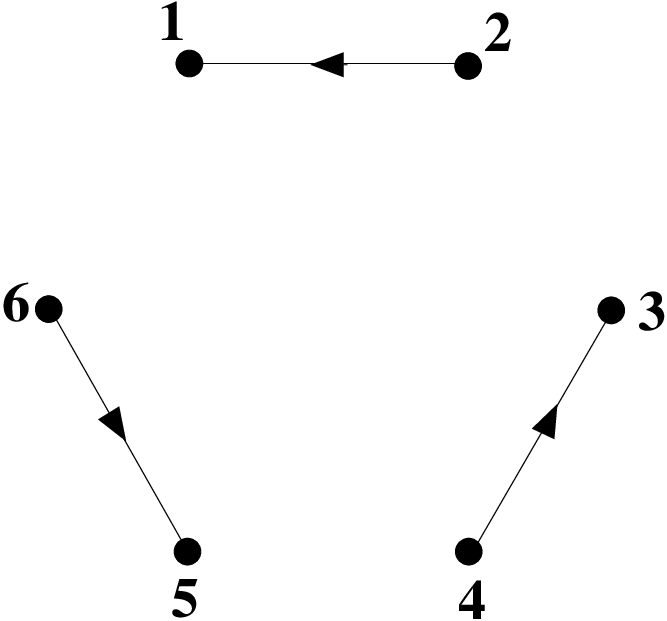}\label{subfig:P6212}
}\qquad
\subfigure[$P_{6,2}^{(2,1)}$]{
\includegraphics[width=3cm]{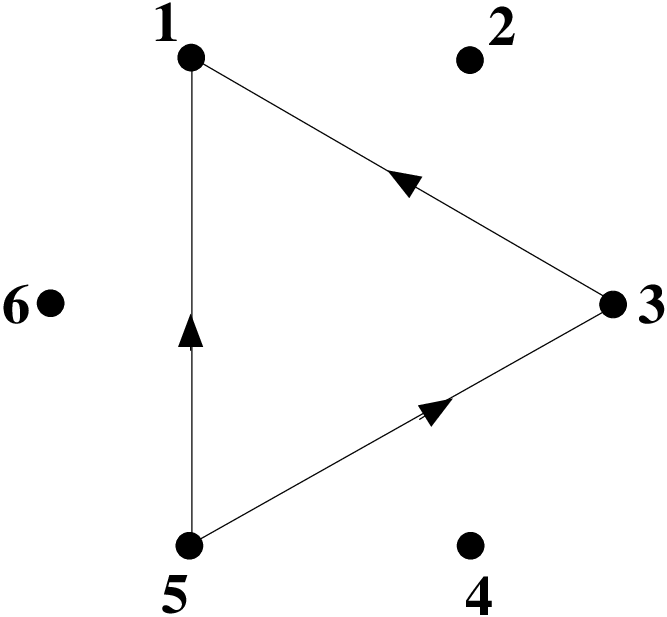}\label{subfig:P6221}
}\qquad
\subfigure[$P_{6,2}^{(2,2)}$]{
\includegraphics[width=3cm]{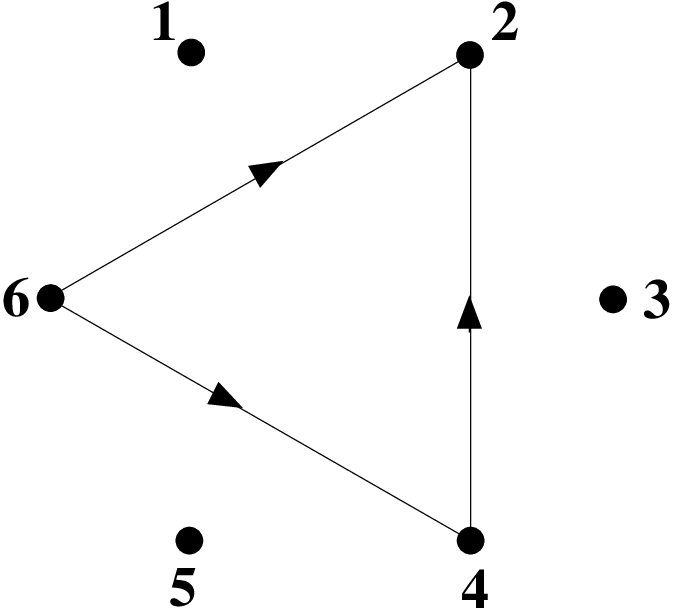}\label{subfig:P6222}
}
\caption{The period $2$ primitives for $6$ nodes}
\label{6node-p2}
\end{figure}

Figures~\ref{4node-p2} and \ref{6node-p2} show the
strictly period $2$ primitives with $4$ and $6$ nodes.

We need the following:

\bl \label{taupreserves} %
\begin{description}
\item[(a)]
Let $M$ be an $N\times N$ skew-symmetric matrix with $M_{ij}\geq 0$ whenever
$i\geq j$. Then $\tau M \tau^{-1}$ has the same property.
\item[(b)]
All period $2$ primitives $B_{N,2}^{(k,l)}$ have nonnegative entries below
the leading diagonal.
\end{description}
\el %

\noindent \textbf{Proof}:
We must also have $M_{ij}\leq 0$ for $i\leq j$. We have
$$(\tau M \tau^{-1})_{ij}=
\begin{cases}
M_{i-1,j-1} & i>1,\ j>1, \\
-M_{N,j-1}  & i=1,\ j>1, \\
-M_{i-1,N}  & i>1,\ j=1, \\
M_{N,N}     & i=1,\ j=1.
\end{cases}$$
from which (a) follows. Part (b) follows from part (a) and the definition
of the period $2$ primitives.
$\Box$

As in the period $1$ case, we obtain period $2$ sink-type quivers by taking
orbit-sums of the basis elements:

\bp[Classification of period $2$ sink-type
quivers]\label{p2sinktypeclassification} If $N$ is odd, there are no strictly
period $2$ sink-type quivers with $N$ nodes. If $N=2r$ is an even integer then
every strictly period $2$ sink-type quiver with $N$ nodes has corresponding
matrix of the form
$$
B =
\left\{ \begin{array}{cc}
\sum_{k=1}^{r}\sum_{j=1}^2 m_{kj} B_{N,2}^{(k,j)} & \mbox{if\ }4|N, \\
(\sum_{k=1}^{r-1}\sum_{j=1}^2 m_{kj} B_{N,2}^{(k,j)})+m_{r1}B_N^{(r)} &
\mbox{if\ }4\nmid N,
\end{array} \right.
$$
where the $m_{jk}$  are arbitrary nonnegative integers such that if $4|N$,
there is at least one $k$, $1\leq k \leq r$, such that $m_{k1}\not=m_{k2}$, and
if $4\nmid N$, there is at least one $k$, $1\leq k\leq r-1$, such that
$m_{k1}\not=m_{k2}$.
\ep  %

\noindent \textbf{Proof}: Using the above discussion and an argument similar
to that in the period $1$ case, we obtain an expression as above for $B$
for which the $m_{kj}$ are integers. It is easy to check that each primitive
has a non-zero entry in the first or second column, below the leading diagonal.
By Lemma~\ref{taupreserves}, this entry must be positive. If the entry is
in the first column, the corresponding $m_{kj}$ must be nonegative as $1$ is
a sink. If it is in the second column then, since $1$ is a sink, mutation
at $1$ does not affect the entries in the second column below the leading
diagonal. Since after mutation at $1$, $2$ is a sink, the corresponding
$m_{kj}$ must be nonnegative in this case also.$\Box$

Whilst the formulae above depend upon particular characteristics of the
primitives, i.e. having a specific sink, a similar relation exists for
\emph{any} period $2$ quiver.
For {\bf any} quiver $Q$
(regardless of any symmetry or periodicity properties), we have
$\mu_{k+1}\,\rho\, Q = \rho\,\mu_k\, Q$, which just corresponds to relabelling
the nodes.  We write this symbolically as $\mu_{k+1}\rho = \rho\mu_k$ and
$\rho^{-1}\mu_{k+1}= \mu_k\rho^{-1}$.  For the period $2$ case, the periodic
chain~(\ref{periodchain}) can be written as:
$$
Q(1) \stackrel{\mu_1}{\longrightarrow} Q(2) \stackrel{\mu_2}{\longrightarrow}
  Q(3)=\rho^2 Q(1)  \stackrel{\mu_3}{\longrightarrow} Q(4)=\rho^2 Q(2)
       \stackrel{\mu_4}{\longrightarrow} \, \cdots
$$
Whilst $\mu_1$ and $\mu_2$ are genuinely different mutations, $\mu_3$ and
$\mu_4$ are just $\mu_1$ and $\mu_2$ after relabelling.
Since $\rho^{-1}\mu_2 Q(2) = \rho Q(1)$,
we have $\mu_1 (\rho^{-1}Q(2))=\rho Q(1)$.

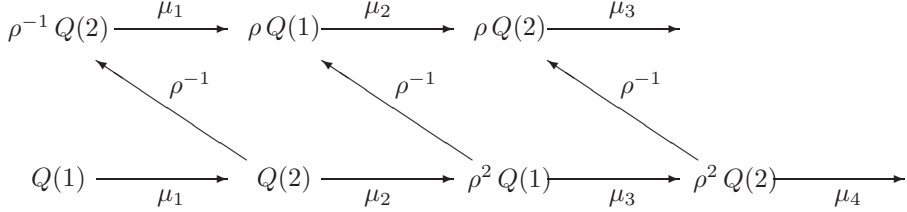
\begin{figure}[hbt]
\begin{center}
\unitlength=0.5mm
\begin{picture}(240,55)
\put(10,45){\makebox(0,0){$\rho^{-1}\, Q(2)$}}
\put(70,45){\makebox(0,0){$\rho\,
Q(1)$}} \put(130,45){\makebox(0,0){$\rho\, Q(2)$}}
\put(10,5){\makebox(0,0){$Q(1)$}} \put(70,5){\makebox(0,0){$Q(2)$}}
\put(130,5){\makebox(0,0){$\rho^2\, Q(1)$}} \put(190,5){\makebox(0,0){$\rho^2\,
Q(2)$}}
\put(25,45){\vector(1,0){30}} \put(80,45){\vector(1,0){35}}
\put(140,45){\vector(1,0){35}} \put(20,5){\vector(1,0){35}}
\put(80,5){\vector(1,0){35}} \put(140,5){\vector(1,0){35}}
\put(200,5){\vector(1,0){35}} \put(60,10){\vector(-3,2){40}}
\put(120,10){\vector(-3,2){40}} \put(180,10){\vector(-3,2){40}}
\put(40,50){\makebox(0,0){$\mu_1$}} \put(95,50){\makebox(0,0){$\mu_2$}}
\put(160,50){\makebox(0,0){$\mu_3$}} \put(40,0){\makebox(0,0){$\mu_1$}}
\put(95,0){\makebox(0,0){$\mu_2$}} \put(160,0){\makebox(0,0){$\mu_3$}}
\put(220,0){\makebox(0,0){$\mu_4$}} \put(45,30){\makebox(0,0){$\rho^{-1}$}}
\put(105,30){\makebox(0,0){$\rho^{-1}$}}
\put(165,30){\makebox(0,0){$\rho^{-1}$}}
\end{picture}
\end{center}
\caption{Period $2$ quivers and mutations} \label{p2quivers}
\end{figure}

We also have $\mu_2\,(\rho\,Q(1)) = \rho\,\mu_1\, Q(1) = \rho\, Q(2)$.
Since $Q(3) = \mu_2 Q(2) = \rho^2 Q(1)$, we have
$\rho^{-1} \mu_2 Q(2) = \rho Q(1)$, and thus we obtain
$\mu_1(\rho^{-1} Q(2)) = \rho Q(1)$.
We thus can extend the above diagram to that in Figure~\ref{p2quivers}.

If $Q(1), Q(2)$ have {\em sinks} at nodes $1$ and $2$ respectively, then so do
$\rho^{-1}Q(2)$ and $\rho Q(1)$ and the mutations $\mu_1$ and $\mu_2$ in the
above diagram act {\em linearly}.  This gives
$$
\mu_1 (Q(1)+\rho^{-1}Q(2)) = Q(2)+\rho\, Q(1) = \rho (Q(1)+\rho^{-1}Q(2))
$$
and
$$\mu_2 (Q(2)+\rho\, Q(1)) = \rho^2\, Q(1) + \rho\, Q(2) = \rho
(Q(2)+\rho\, Q(1)),
$$
so $Q(1)+\rho^{-1}Q(2)$ is period $1$.

We have proved the following:

\bp  %
Let $Q$ be period $2$ sink-type quiver. Then $Q(1)+\rho^{-1}Q(2)$ is a quiver
of period $1$.
\ep  %

\section{Quivers with Higher Period}\label{p-high}

Higher period primitives are defined in a similar way.  The
periodic chain (\ref{periodchain}) contains $m$ quivers
$Q(1),Q(2),\ldots ,Q(m)$, with corresponding matrices $B(1),\cdots ,B(m)$.

\bd[Period $m$ sink-type quivers]\label{d:periodm}   %
A quiver $Q$ is said to be a \emph{period $m$ sink-type quiver} if it is of
period $m$ and, for $1\leq i\leq m$, node $i$ of $Q(i)$ is a sink.
\ed   %

Thus the mutation $Q(i)\mapsto Q(i+1)=\mu_i Q(i)$ again only involves the
reversal of arrows at node $i$, so can be achieved through a simple conjugation
of its matrix: $\mu_i B(i) = D_i B(i) D_i$. Here
$$
D_i= \mbox{diag}(1,\cdots ,1,-1,1,\cdots ,1) = \rho^{i-1} D_1
\rho^{-i+1}
$$
(with a ``$-1$'' in the $i$th position).

As in the period $1$ and $2$ cases, we obtain:

\bl[Period $m$ sink-type equation]
Suppose that $Q$ is a quiver with a sink at the $i$th node of $Q(i)$ for
$i=1,2,\ldots ,m$. Then $Q$ is period $m$ if and only if
$\tau^m B_Q \tau^{-m}=B_Q$.
\el

\noindent \textbf{Proof}:
We have that $Q$ has period $m$ if and only if
$D_m\cdots D_1 B_Q D_1\cdots D_m=\rho^m B_Q \rho^{-m}$, i.e.\
if and only if
$$
B_Q = D_1\cdots D_m \rho^m B_Q \rho^{-m}D_m\cdots D_1 =
         \tau^m B_Q \tau^{-m}.
$$
$\Box$

Starting with the same matrices $R_N^{(k)}$, we now use the action $M\mapsto
\tau^m M\tau^{-m}$ to build an invariant, i.e.\ we take orbit sums for
$\tau^m$. We only obtain strictly $m$-periodic elements in the case where the
orbit has size divisible by $m$.

When $m|N$, the matrices $R_N^{(k)}$, for $1\leq k\leq r-1$ (where $N=2r$ or
$2r+1$, $r$ an integer), generate period $m$ primitives $B_{N,m}^{(k,1)}$, with
matrices given by
$$
B_{N,m}^{(k,1)} = \sum_{i=0}^{(N/m)-1} \tau^{mi} R_N^{(k)} \tau^{-mi}.
$$
Geometrically, the primitive $P_{N,m}^{(k,1)}$ is obtained from the primitive
$P_N^{(k)}$ by only including every $m^{th}$ arrow.  As before, we form another
$m-1$ period $m$ primitives, $P_{N,m}^{(k,j)}$ for $j=2,\ldots ,m$, with
matrices given by:
$$
B_{N,m}^{(k,j)} =  \tau^{j-1} B_{N,m}^{(k,1)} (\tau^{j-1})^{-1}.
$$
Note that the elements $\tau^l R_N^{(k)} \tau^{-l}$,
for $0\leq l\leq N-1$, form a $\tau$-orbit of size $N$. Since $m|N$, this
breaks up into $m$ $\tau^m$-orbits each of size $N/m$; the elements above
are the orbit sums.

Similarly, if $(2m)|N$ (so we are in the case $N=2r$) then the
$\tau^m$-orbit-sum of $R_N^{(r)}$ is
$$
B_{N,m}^{(k,1)} = \sum_{i=0}^{(N/2m)-1} \tau^{mi} R_N^{(r)} \tau^{-mi}.
$$
with corresponding quiver $P_{N,m}^{(k,1)}$. We also obtain another $m-1$
period $m$ primitives, $P_{N,m}^{(r,j)}$, for $j=2,\ldots ,m$, with matrices
$$
B_{N,m}^{(r,j)} =  \tau^{j-1} B_{N,m}^{(r,1)} (\tau^{j-1})^{-1}.
$$

As in the period $1$ and $2$ cases, we obtain arbitrary strictly period $m$
sink-type quivers by taking orbit-sums of the basis elements. The nonnegativity
of the coefficients $m_{kj}$ is shown in a similar way also.

\bp[Classification of period $m$ sink-type
quivers]\label{pmsinktypeclassification} If $m\nmid N$, there are no strictly
period $m$ sink-type quivers. If $(2m)|N$, the general strictly period $m$
sink-type quiver is of the form
$$B = \sum_{k=1}^{r}\sum_{j=1}^{m} m_{kj} B_{N,m}^{(k,j)},$$
where the $m_{kj}$ are nonnegative integers and there is at least one $k$,
$1\leq k\leq r$, for which the $m_{kj}$ are not all equal.

If $m|N$ but $(2m)\nmid N$ then the general period $m$ sink-type quiver has the
form
$$B =
\left\{ \begin{array}{cc} \sum_{k=1}^{r}\sum_{j=1}^{m} m_{kj} B_{N,m}^{(k,j)} &
\textrm{ if $N=2r+1$ is
odd;} \\
\sum_{k=1}^{r-1}\sum_{j=1}^{m} m_{kj} B_{N,m}^{(k,j)}+ \sum_{j=1}^{m/2} m_{rj}
B_{N,m/2}^{(r,j)}
  & \textrm{ if $N=2r$ is even,}
\end{array} \right.
$$
where the $m_{kj}$ are nonnegative integers and where in the first case,
there is at least one $k$, $1\leq k\leq r$, for which
the $m_{kj}$ are not all equal, and in the second case, there is at least one
$k$, $1\leq k\leq r-1$, for which the $m_{kj}$ are not all equal.
\ep  %

As before, we use $\mu_{k+1}\rho = \rho\mu_k$ and $\rho^{-1}\mu_{k+1}=
\mu_k\rho^{-1}$, from which it follows that
$\mu_k\rho^{-j}=\rho^{-j}\mu_{k+j}$.  In turn, this gives
$$
\mu_k\, (\rho^{-j}\, Q(j+k)) = \rho^{-j}\,\mu_{j+k}\, Q(j+k)
                = \rho^{-j}\, Q(j+k+1).
$$
Suppose now that $Q$ is a period $m$ quiver.
Then we have $Q(sm+j)=\rho^{sm}\, Q(j)$ for $1\leq j\leq m$. We use this to
extend the periodic chain~(\ref{periodchain}) to an $m$ level array.
We have
$$
\mu_1\, (\rho^{-j}\, Q(j+1))=\rho^{-j}\, Q(j+2),\;\; \mu_2\, (\rho^{-j}\,
Q(j+2))=\rho^{-j}\, Q(j+3), \ldots ,
$$
arriving at
$$
\mu_m\, (\rho^{-j}\, Q(j+m) =\rho^{-j}\, Q(j+m+1)=\rho^m (\rho^{-j}\,
Q(j+1)).
$$
We write this period $m$ sequence in the $j$th level of the array, i.e.\
$$
\rho^{-j}\, Q(j+1) \stackrel{\mu_1}{\longrightarrow} \rho^{-j}\, Q(j+2)
\stackrel{\mu_2}{\longrightarrow}
  \cdots \stackrel{\mu_{m-1}}{\longrightarrow} \rho^{-j}\, Q(j+m)
       \stackrel{\mu_m}{\longrightarrow} \rho^m (\rho^{-j}\,
Q(j+1)) .
$$

Again we have that if $Q(j)$ has a sink at node $j$ for each $j$, then
each $\rho^{-j}Q(j+1)$ has a sink at node $1$ and the mutation $\mu_1$ acts
linearly. This gives
$$
\mu_1 (Q(1)+\rho^{-1}Q(2)+\cdots +\rho^{-m+1}Q(m)) = \rho
(Q(1)+\rho^{-1}Q(2)+\cdots +\rho^{-m+1}Q(m)),
$$
so $Q(1)+\rho^{-1}Q(2)+\cdots +\rho^{-m+1}Q(m)$ is period $1$.

We have proved:

\bp
Let $Q$ be period $m$ sink-type quiver.
Then $Q(1)+\rho^{-1}Q(2)+\cdots +\rho^{-m+1}Q(m)$ is a quiver of period $1$.
\ep

\begin{figure}[htb]
\centering
\subfigure[$P_{6,3}^{(1,1)}$]{
\includegraphics[width=3cm]{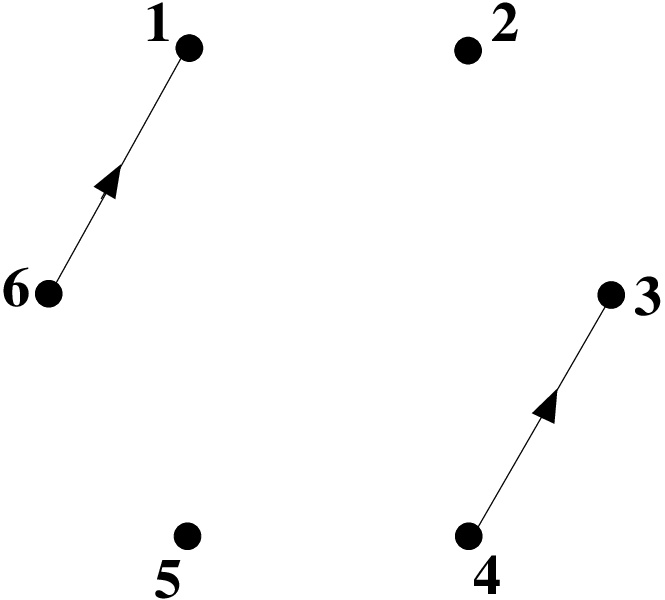}\label{subfig:P6311}
}\qquad
\subfigure[$P_{6,3}^{(1,2)}$]{
\includegraphics[width=3cm]{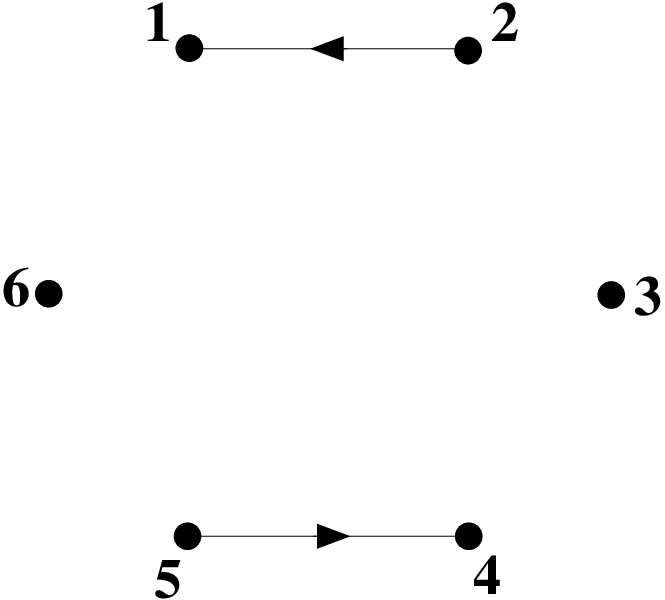}\label{subfig:P6312}
}\qquad
\subfigure[$P_{6,3}^{(1,3)}$]{
\includegraphics[width=3cm]{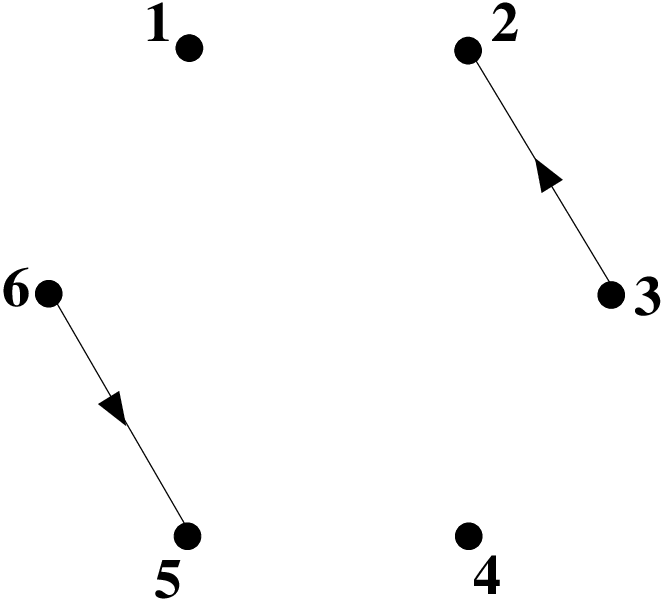}\label{subfig:P6313}
}\\
\subfigure[$P_{6,3}^{(2,1)}$]{
\includegraphics[width=3cm]{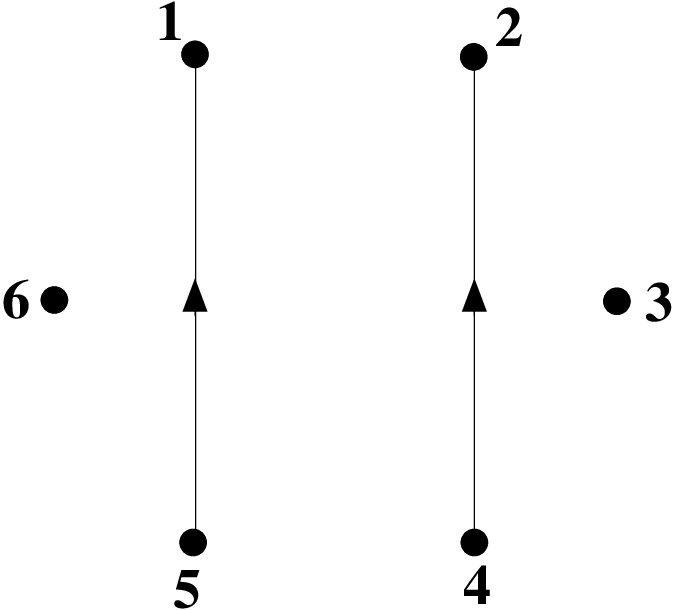}\label{subfig:P6321}
}\qquad
\subfigure[$P_{6,3}^{(2,2)}$]{
\includegraphics[width=3cm]{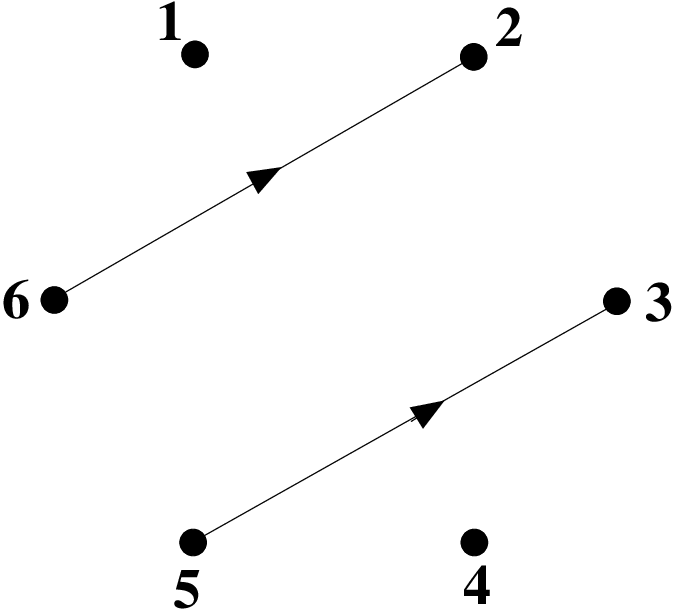}\label{subfig:P6322}
}\qquad
\subfigure[$P_{6,3}^{(2,3)}$]{
\includegraphics[width=3cm]{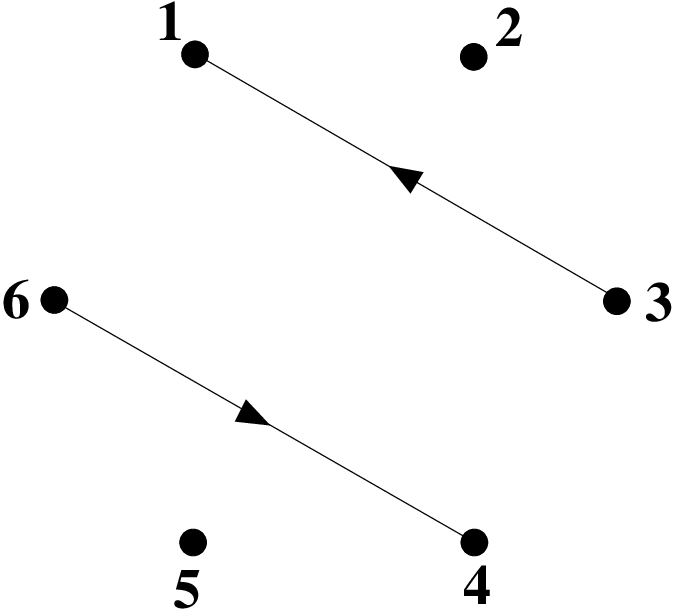}\label{subfig:P6323}
}\\
\subfigure[$P_{6,3}^{(3,1)}$]{
\includegraphics[width=3cm]{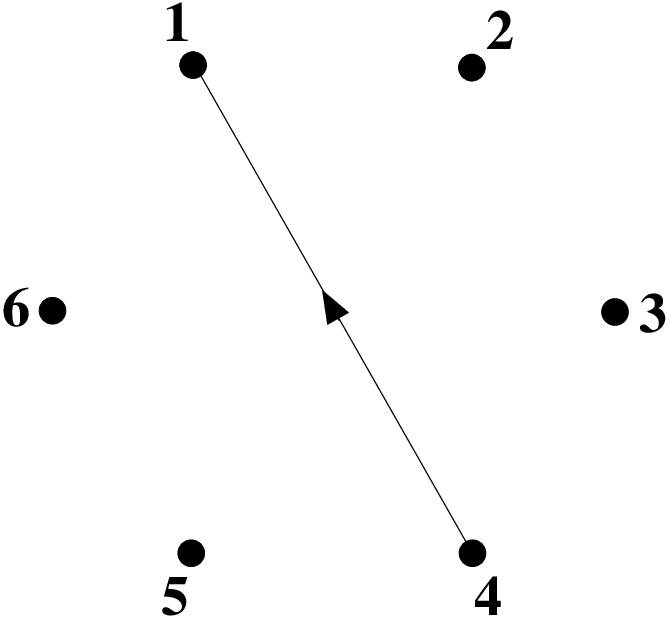}\label{subfig:P6331}
}\qquad
\subfigure[$P_{6,3}^{(3,2)}$]{
\includegraphics[width=3cm]{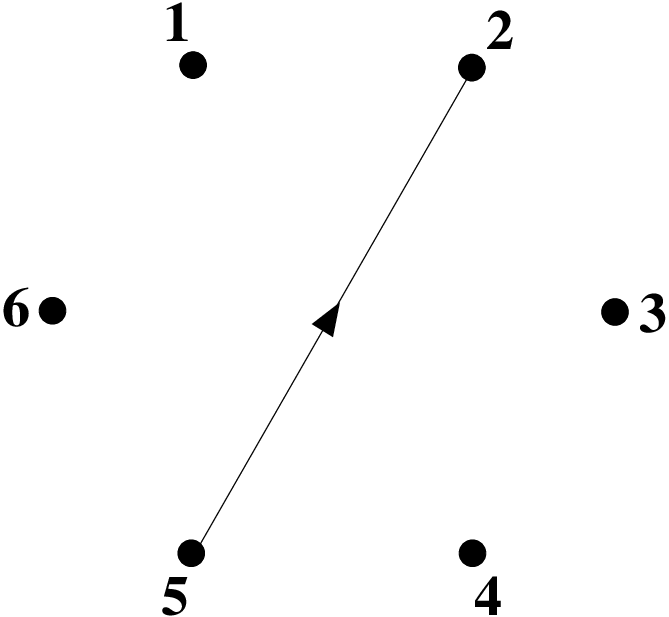}\label{subfig:P6332}
}\qquad
\subfigure[$P_{6,3}^{(3,3)}$]{
\includegraphics[width=3cm]{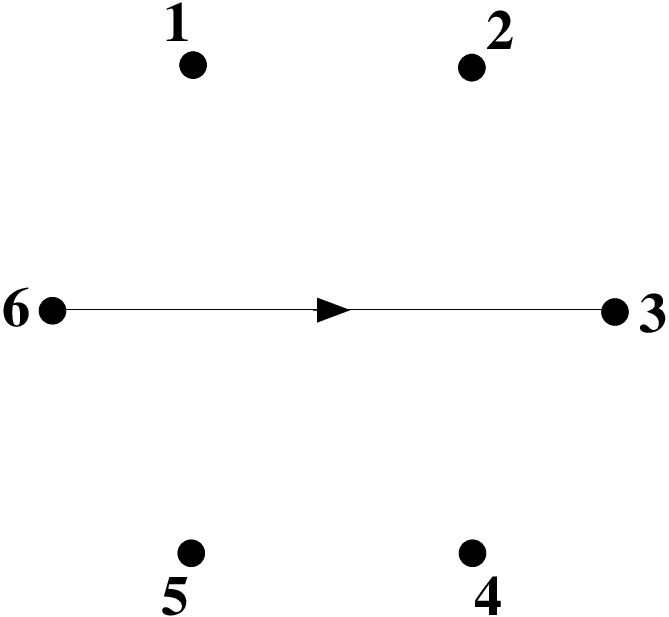}\label{subfig:P6333}
}
\caption{The period $3$ primitives for $6$ nodes.}
\label{6node-p3}
\end{figure}

\bex[Period $3$ Primitives]  \label{p3primitives}  {\em   %
Proceeding as described above, whenever $N$ is a multiple of $3$ we obtain $3$
period $3$ primitives for each period $1$ primitive. Figure \ref{6node-p3}
shows those with $6$ nodes.
}\eex   %

\section{Period $1$ General Solution}\label{p1-gen}

In this section we give an explicit construction of the $N\times N$
skew-symmetric matrices corresponding to arbitrary period $1$ quivers,
i.e.\ those for which mutation at node $1$ has the same effect as the
rotation $\rho$. We express the general solution as an explicit sum
of period $1$ primitives, thus giving a simple classification of all such
quivers.

In anticipation of the final result, we consider the following matrix:
\be  \label{genp1-matrix}    %
B = \left(\begin{array}{ccccc}
     0 & -m_1 &  \cdots & -m_{N-1} \\
     m_1 & 0  && *  \\
     \vdots && 0 &  \\
     m_{N-1} & * && 0
     \end{array}  \right).
\ee    %
Using (\ref{gen-mut}), the general mutation rule at node $1$ is
\be   \label{mut1}  %
\tilde b_{ij}= \left\{ \begin{array}{ll}
                        -b_{ij} & \mbox{if}\;\; i=1\;\;\mbox{or}\;\; j=1, \\
                        b_{ij}+\frac{1}{2} (|b_{i1}|b_{1j}+b_{i1}|b_{1j}|)
                        & otherwise.
                        \end{array}  \right.
\ee    %
The effect of the rotation $B\mapsto \rho B \rho^{-1}$ is to move the
entries of $B$ down and right one step, so that
$(\rho B \rho^{-1})_{ij}=b_{i-1,j-1}$, remembering that indices are
labelled modulo $N$, so $N+1\equiv 1$.
For $1\leq i,j\leq N-1$, let
$$\varepsilon_{ij}=\frac{1}{2}(m_i|m_j|-m_j|m_i|).$$
Then if $m_i$ and $m_j$ have the same sign, $\varepsilon_{ij}=0$. Otherwise
$\varepsilon_{ij}=\pm |m_im_j|$, where the sign is that of $m_i$. Let
$\widetilde{B}=\mu_1 B$, so that $\tilde{b}_{ij}=b_{ij}+\varepsilon_{i-1,j-1}$.

\bt \label{p1-prototheorem}    %
Let $B$ be an $N\times N$ skew-symmetric integer matrix. Let $b_{k1}=m_{k-1}$
for $k=2,3,\ldots ,N$. Then $\mu_1 B=\rho B\rho^{-1}$ if and only if
$m_r=m_{N-r}$ for $r=1,2,\ldots ,N-1$,
$b_{ij}=m_{i-j}+\varepsilon_{1,i-j+1}+\varepsilon_{2,i-j+2}+\cdots +\varepsilon_{j-1,i-1}$ for
all $i>j$, and $B$ is symmetric along the non-leading diagonal.
\et     %

\noindent \textbf{Proof}:
By skew-symmetry, we note that we only need to determine $b_{ij}$ for $i>j$.
We need to solve $\mu_1 B=\rho B\rho^{-1}$.
By the above discussion, this is equivalent to solving
\begin{equation} \label{e:mutationequal}
b_{ij}+\varepsilon_{i-1,j-1} = b_{i-1,j-1},
\end{equation}
for $i>j$, with $\varepsilon_{ij}$ as given above. Solving the equation leads to a
recursive formula for $b_{ij}$.

We obtain
\begin{eqnarray*} \label{e:bij}
b_{ij} &=& b_{i-1,j-1}+\varepsilon_{j-1,i-1} \\
&=& b_{i-2,j-2}+\varepsilon_{j-1,i-1}+\varepsilon_{j-2,i-2} \\
&\vdots & \\
&=& b_{i-j+1,1}+\varepsilon_{j-1,i-1}+\varepsilon_{j-2,i-2}+\cdots +\varepsilon_{1,i-j+1}.
\end{eqnarray*}

In particular, we have:
\begin{equation} \label{e:bnj}
b_{Nj} = m_{N-j}+\varepsilon_{1,N-j+1}+\varepsilon_{2,N-j+2}+\cdots +
     \varepsilon_{j-2,N-2}+\varepsilon_{j-1,N-1}.
\end{equation}
We also have that $m_j=\tilde{b}_{1,j+1}=(\rho B \rho^{-1})_{1,j+1}=b_{Nj}$.
In particular, $m_1=b_{N1}=m_{N-1}$.
Equation~\ref{e:bnj} gives
$$
m_2=b_{N2}=m_{N-2}+\varepsilon_{1,N-1}=m_{N-2}+\varepsilon_{11}=m_{N-2}.
$$
So $m_2=m_{N-2}$.
Suppose that we have shown that $m_j=m_{N-j}$ for $j=1,2,\ldots ,r$.
Then equation~\ref{e:bnj} gives
\begin{eqnarray*}
b_{N,r+1} &=& m_{N-r-1}+\varepsilon_{1,N-r}+\varepsilon_{2,N-r+1}+\cdots +\varepsilon_{r,N-1} \\
&=& m_{N-r-1}+\sum_{i=1}^r \varepsilon_{i,N-r+i-1} \\
&=& m_{N-r-1}+\sum_{i=1}^r \varepsilon_{i,r+1-i} \\
&=& m_{N-r-1}+\varepsilon_{1,r}+\varepsilon_{2,r-1}+\cdots +\varepsilon_{r,1} = m_{N-r-1},
\end{eqnarray*}
using the inductive hypothesis and the fact that $\varepsilon_{st}=-\varepsilon_{ts}$
for all $s,t$.
Hence $m_{r+1}=m_{N-r-1}$ and we have by induction that $m_r=m_{N-r}$ for
$1\leq r\leq N-1$.

We have, for $i>j$, by equation~(\ref{e:bij}),
\begin{eqnarray*}
b_{N-j+1,N-i+1} &=& m_{(N-j+1)-(N-i+1)}+\varepsilon_{(N-i+1)-1,(N-j+1)-1}+
   \varepsilon_{(N-i+1)-2,(N-j+1)-2}+\cdots +\varepsilon_{1,(N-j+1)-(N-i+1)+1} \\
&=& m_{i-j}+\varepsilon_{N-i,N-j}+\varepsilon_{N-i-1,N-j-1}+\cdots +\varepsilon_{1,i-j+1},
\end{eqnarray*}
and we have, again using~(\ref{e:bij}) and the fact that
$\varepsilon_{N-a,N-b}=\varepsilon_{ab}$,
\begin{eqnarray*}
m_{i-j}=m_{N-i+j} &=& b_{N,N-i+j} \\
&=& b_{N-j+1,N-i+1}+\varepsilon_{N-i+j-1,N-1}+\varepsilon_{N-i+j-2,N-2}+\cdots
+\varepsilon_{N-i+1,N-j+1} \\
&=& b_{N-j+1,N-i+1}+\varepsilon_{i-j+1,1}+\varepsilon_{i-j+2,2}+\cdots
+\varepsilon_{i-1,j-1},
\end{eqnarray*}
so
$$b_{N-j+1,N-i+1}=m_{i-j}+\varepsilon_{j-1,i-1}+\cdots +\varepsilon_{i-j+2,2}+
\varepsilon_{1,i-j+1}=b_{ij}.$$
Hence $B$ is symmetric along the non-leading diagonal.

If $B$ satisfies all the requirements in the statement of the theorem,
then equation~\ref{e:bij} is satisfied, and therefore
$\rho B \rho^{-1}=\mu_1 B$. The proof is complete.$\Box$

We remark that with the identification $m_r=m_{N-r}$, we have seen that
the formula~(\ref{e:bnj}) has a symmetry, due to which the $\varepsilon$'s cancel
in pairwise fashion:
$$
b_{N,N-k+1} = b_{k1}+\varepsilon_{1k}+\varepsilon_{2,k+1}+\cdots +
     \varepsilon_{k+1,2}+\varepsilon_{k1}
$$
The formula (\ref{e:bij}) is just a truncation of this, so not all terms
cancel.  As we march from $b_{k1}$ in a ``south easterly direction'', we first
add $\varepsilon_{1k}, \varepsilon_{2,k+1}$, etc, until we reach $\varepsilon_{r,r+1}$
(when $N-k=2r$) or $\varepsilon_{rr}=0$ (when $N-k=2r+1$).  At this stage we start
to subtract terms on a basis of ``last in, first out'', with the result that
the matrix has reflective symmetry about the second diagonal as we have
seen.

\br[Sink-type case]
We note that if all the $m_i$ have the same sign, then all the
$\varepsilon_{ij}$ are zero. Equation~(\ref{e:mutationequal}) reduces to
$b_{ij}=b_{i-1,j-1}$ and we recover the sink-type period $1$ solutions
considered in Proposition~\ref{p1sinktypeclassification}.
\er

\subsection{Examples}  \label{p1-examples}

The simplest nontrivial example is when $N=4$.

\bex[Period $1$ Quiver with $4$ Nodes]  \label{n=4p1}  {\em %
Here the matrix has the form
$$
B = \left(\begin{array}{cccc}
     0 & -m_1 & -m_2 &  -m_1 \\
     m_1 & 0 & -m_1-\varepsilon_{12} & -m_2 \\
     m_2 & m_1+\varepsilon_{12} & 0 & -m_1 \\
     m_1 & m_2 & m_1 & 0
     \end{array}\right) ,
$$

As previously noted, if $m_1$ and $m_2$ have the same sign, then
$\varepsilon_{12}=0$ and this matrix is just the sum of primitives for $4$ nodes.
The $2\times 2$ matrix in the ``centre'' of $B$ (formed out of rows and
columns $2$ and $3$),
$$
  \left(\begin{array}{cc}
     0 & -\varepsilon_{12} \\
     \varepsilon_{12} & 0
     \end{array}\right),
$$
corresponds to $\varepsilon_{12}$ times the primitive $P_2^{(1)}$ with $2$ nodes
(see Figure \ref{234node}). For the case $m_1=1, m_2=-2, \varepsilon_{12}=2$, we
obtain the Somos $4$ quiver in Figure~\ref{subfig:somos4quiver}. The action of
$\iota$ (see Remark~\ref{iota}) is $1\leftrightarrow 4, 2\leftrightarrow 3$ and
clearly just reverses all the arrows as predicted by Remark~\ref{iota}. }\eex

\bex[Period $1$ Quiver with $5$ Nodes]  \label{n=5p1}  {\em %
Here the general period $1$ solution has the form
$$
B = \left(\begin{array}{ccccc}
     0 & -m_1 & -m_2 & -m_2 & -m_1 \\
     m_1 & 0 & -m_1-\varepsilon_{12} & -m_2-\varepsilon_{12} & -m_2 \\
     m_2 & m_1+\varepsilon_{12}  & 0 & -m_1-\varepsilon_{12}  & -m_2 \\
     m_2 & m_2+\varepsilon_{12}  & m_1+\varepsilon_{12} & 0 & -m_1 \\
     m_1 & m_2 & m_2 & m_1 & 0
     \end{array}\right)
$$
which can be written as
$$
B = \sum_{k=1}^2 m_k\, B_5^{(k)} +  \varepsilon_{12} B_3^{(1)} ,
$$
where $B_3^{(1)}$ is embedded symmetrically in the middle of a
$5\times 5$ matrix (surrounded by zeros).

When $m_1=1$ and $m_2=-1$, this matrix corresponds to the Somos $5$ sequence;
see Figure~\ref{fig:somos5quiver} for the corresponding quiver.
\begin{figure}[ht]
\centering
\includegraphics[width=4cm]{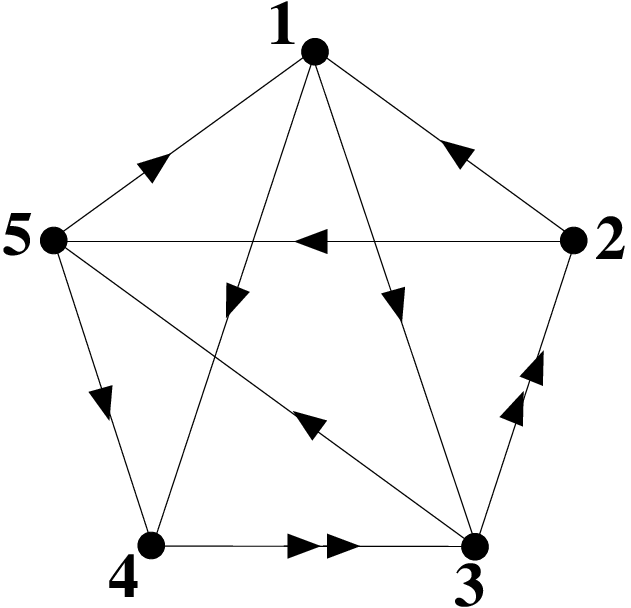}
\caption{The Somos $5$ quiver.}\label{fig:somos5quiver}
\end{figure}
}\eex   %

\bex[Period $1$ Quiver with $6$ Nodes]  \label{n=6p1}  {\em %
Here the matrix has the form
\bea   %
B &=& \left(\begin{array}{cccccc}
     0 & -m_1 & -m_2 & -m_3 & -m_2 & -m_1 \\
     m_1 & 0 & -m_1 & -m_2 & -m_3 & -m_2 \\
     m_2 & m_1 & 0 & -m_1 & -m_2 & -m_3 \\
     m_3 & m_2 & m_1 & 0 & -m_1 & -m_2 \\
     m_2 & m_3 & m_2 & m_1 & 0 & -m_1 \\
     m_1 & m_2 & m_3 & m_2 & m_1 & 0
     \end{array}\right)   \nn\\[3mm]
     && \quad +
     \left(\begin{array}{c|cccc|c}
     0 & 0 & 0 & 0 & 0 & 0 \\
     \hline
     0 & 0 & -\varepsilon_{12} & -\varepsilon_{13} & -\varepsilon_{12} & 0 \\
     0 & \varepsilon_{12} & 0 & -\varepsilon_{12} & -\varepsilon_{13} & 0 \\
     0 & \varepsilon_{13} & \varepsilon_{12} & 0 & -\varepsilon_{12} & 0 \\
     0 & \varepsilon_{12} & \varepsilon_{13} & \varepsilon_{12} & 0 & 0 \\
     \hline
     0 & 0 & 0 & 0 & 0 & 0
     \end{array}\right)  +
     \left(\begin{array}{cc|cc|cc}
     0 & 0 & 0 & 0 & 0 & 0 \\
     0 & 0 & 0 & 0 & 0 & 0 \\
     \hline
     0 & 0 & 0 & -\varepsilon_{23} & 0 & 0 \\
     0 & 0 & \varepsilon_{23} & 0 & 0 & 0 \\
     \hline
     0 & 0 & 0 & 0 & 0 & 0 \\
     0 & 0 & 0 & 0 & 0 & 0
     \end{array}\right)   , \nn
\eea   %
which can be written as
$$
B = \sum_{j=1}^3 m_k\, B_6^{(k)} + \sum_{k=1}^2 \varepsilon_{1,k+1}\, B_4^{(k)} +
   \varepsilon_{23}\, B_2^{(1)} ,
$$
where the periodic solutions with fewer rows and columns are embedded
symmetrically within a $6\times 6$ matrix.
}\eex   %

\subsection{The Period $1$ General Solution in Terms of Primitives}

It can be seen from the above examples that the solutions are built out of a
sequence of sub-matrices, each of which corresponds to one of the primitives.
The main matrix is just an integer linear combination of primitive matrices
for the full set of $N$ nodes.  The next
matrix is a combination (with coefficients $\varepsilon_{1j}$) of primitive matrices
for the $N-2$ nodes $2, \cdots, N-1$.  We continue to reduce by $2$ until we
reach either $2$ nodes (when $N$ is even) or $3$ nodes (when $N$ is odd).

Remarkably, as can be seen from the general structure of the matrix given by
(\ref{e:bij}), together with the symmetry $m_{N-r}=m_r$, this description
holds for all $N$.

Recall that for an even (or odd) number of nodes, $N=2r$ (or
$N=2r+1$), there are $r$ {\em primitives}, labelled $B_{2r}^{(k)}$ (or
$B_{2r+1}^{(k)}$), $k=1,\cdots ,r$.
We denote the general linear combination of these by
$$
\widetilde{B}_{2r}(\mu_1,\cdots ,\mu_r) = \sum_{j=1}^r \mu_j B_{2r}^{(j)},
\quad\mbox{or}\quad
\widetilde{B}_{2r+1}(\mu_1,\cdots ,\mu_r) = \sum_{j=1}^r \mu_j
B_{2r+1}^{(j)},
$$
for integers $\mu_j$.

The quivers corresponding to
$\widetilde{B}_{2r}(\mu_1,\cdots ,\mu_r)$ and
$\widetilde{B}_{2r+1}(\mu_1,\cdots ,\mu_r)$
(i.e.\ without the extra terms coming from the $\varepsilon_{ij}$)
do not have periodicity properties (in general).

We now restate Theorem~\ref{p1-prototheorem} in this new notation:

\bt[The general period $1$ quiver]  \label{p1-theorem}
Let $B_{2r}$ (respectively $B_{2r+1}$) denote the matrix corresponding
to the general even (respectively odd) node quiver of mutation periodicity
$1$.  Then
\begin{enumerate}
\item
$$
B_{2r} = \widetilde{B}_{2r}(m_1,\cdots ,m_r) +
  \sum_{k=1}^{r-1} \widetilde{B}_{2(r-k)}(\varepsilon_{k,k+1},\cdots ,\varepsilon_{kr}) ,
$$
where the matrix
$\widetilde{B}_{2(r-k)}(\varepsilon_{k,k+1},\cdots ,\varepsilon_{kr})$ is embedded
in a $2r\times 2r$ matrix in rows and columns $k+1,\cdots ,2r-k$.    %
\item
$$
B_{2r+1} = \widetilde{B}_{2r+1}(m_1,\cdots ,m_r) +
  \sum_{k=1}^{r -1}
\widetilde{B}_{2(r-k)+1}(\varepsilon_{k,k+1},\cdots ,\varepsilon_{kr}) ,
$$
where the matrix
$\widetilde{B}_{2(r-k)+1}(\varepsilon_{k,k+1},\cdots ,\varepsilon_{kr})$ is
embedded in a $2r+1\times 2r+1$ matrix in
rows and columns $k+1,\cdots ,2r+1-k$.    %
\end{enumerate}
\et  %

\section{Quivers with Mutation Periodicity $2$}\label{p2-quivers}

Already at period $2$, we cannot give a full classification of the possible
quivers. However, we can give the {\em full list} for low values of $N$, the
number of nodes.  We can also give a class of period $2$ quivers which exists
for odd or even $N$.

When $N$ is {\em even}, primitives play a role, but the full solution cannot be
written purely in terms of primitives.  When $N$ is {\em odd}, primitives do
not even exist, but there are still quivers with mutation periodicity $2$.

Consider the period $2$ chain:
$$
Q(1) \stackrel{\mu_1}{\longrightarrow} Q(2) \stackrel{\mu_2}{\longrightarrow}
   Q(3)=\rho^2 Q(1)
$$
A simpler way to compute is to use $\mu_2 Q(3) = Q(2)$, so
$\mu_2 \rho^2 Q(1) = Q(2)$. Hence we must solve
\be   \label{p2-eqs}  %
    \rho \mu_1 \rho Q(1) = \mu_1 Q(1) ,
\ee   %
which are the equations referred to below.  We first consider the solution of
these equations for $N=3,\cdots , 5$.

We need one new piece of notation, which generalises our former $\varepsilon_{ij}$.
We define
$$
\varepsilon(x,y) = \frac{1}{2} \, (x |y|-y |x|).
$$
Thus, $\varepsilon_{ij}=\varepsilon(m_i,m_j)$.

\subsection{$3$ Node Quivers of Period $2$}

Let

$$B(1) = \left(\begin{array}{ccc}
0 & -m_1 & -m_2 \\
m_1 & 0 & -b_{32} \\
m_2 & b_{32} & 0
\end{array}\right)
$$
Equation (9) gives the equalities $m_1=b_{32}$ and $m_2 - m_1 = -\varepsilon_{12}$.
If the signs of $m_1$ and $m_2$ are the same, we obtain a period $1$ solution.
Assuming otherwise leads to the equation $m_2 - m_1 = \pm m_1 m_2$ depending on the
sign of $m_1$ (and $m_2$). The only integer solutions to this equation are
$m_1 = \pm 2$ and $m_2 = \mp 2$.

It follows that there are just two solutions of period two: the following
matrix and its negative:
$$
B(1) = \left(\begin{array}{ccc}
     0 & -2 & 2 \\
     2 & 0 & -2 \\
     -2 & 2 & 0
     \end{array}\right) .
$$
This corresponds to a $3-$cycle of double arrows.  Notice that in this case,
there are no free parameters.  Mutating at node $1$ just gives $B(2)=-B(1)$,
i.e.\ $Q(2)=Q(1)^{opp}$. Note that the representation theoretic properties of
this quiver are discussed at some length in~\cite[\S 8,\S 11]{08-1}.

\subsection{$4$ Node Quivers of Period 2}\label{p2-4node}

We start with the matrix
$$
B(1) = \left(\begin{array}{cccc}
     0 & -m_1 & -m_2 & -m_3 \\
     m_1 & 0 & -b_{32} & b_{42} \\
     m_2 & b_{32} & 0 & -b_{43} \\
     m_3 & b_{42} & b_{43} & 0
     \end{array}\right)
$$
Setting $p_1=b_{42}$ and solving for $b_{43}$ and $b_{32}$ in terms of $p_1$ and the
$m_i$'s we find:
$$
b_{43}=m_1 ,\quad\mbox{and}\quad b_{32}=m_3 + \varepsilon(m_1,p_1).
$$
We also obtain the $3$ conditions:
$$
\varepsilon_{13}=0,\quad \varepsilon_{12}-\varepsilon(m_1,p_1)=0,\quad
               \varepsilon_{23}+\varepsilon(m_3,p_1)=0.
$$
The first of these three conditions just means that $m_1$ and $m_3$ have the same sign (or that
one of them is zero).  Choosing $m_1>0$, so $m_3\geq 0$, we must have $m_2<0$
for node $1$ {\em not} to be a sink.  The remaining conditions are then
$$
m_1(|p_1|-p_1+2 m_2)=0 ,\quad  m_3(|p_1|-p_1+2 m_2)=0 .
$$
For a nontrivial solution we must have $p_1<0$, which leads to $p_1=m_2$.  The
final result is then the following:
\bea     %
B(1) &=& \left(\begin{array}{cccc}
     0 & -m_1 & -m_2 & -m_3 \\
     m_1 & 0 & m_1m_2-m_3 & -m_2 \\
     m_2 & m_3-m_1m_2 & 0 & -m_1 \\
     m_3 & m_2 & m_1 & 0
     \end{array}\right) ,   \nn\\
     && \label{p2n4}   \\  %
B(2) &=& \left(\begin{array}{cccc}
     0 & m_1 & m_2 & m_3 \\
     -m_1 & 0 & -m_3 & -m_2 \\
     -m_2 & m_3 & 0 & m_2m_3-m_1 \\
     -m_3 & m_2 & m_1-m_2m_3 & 0
     \end{array}\right) ,  \nn
\eea    %
with $m_1>0$, $m_2<0$ and $m_3\geq 0$.
Notice that $B(2)(m_1,m_3)=\rho B(1)(m_3,m_1) \rho^{-1}$, so the period $2$
property stems from the involution $m_1\leftrightarrow m_3$.  If $m_3=m_1$,
then the quiver has mutation period $1$.  We may choose either of these to be
zero, but {\em not} $m_2$, since, again, node $1$ would be a sink.

\br[The Quiver and its Opposite]  \label{qoppREM} %
We made the choice that $m_1>0$.  The equivalent choice $m_1<0$ would just lead
to the negative of $B(1)$, corresponding to $Q(1)^{opp}$.
\er   %

\br[A Graph Symmetry]  \label{graphsymmREM} %
Notice that all $4$-node quivers of period $2$ have the {\bf graph symmetry}
$(1,2,3,4)\leftrightarrow (4,3,2,1)$, under which $Q\mapsto Q^{opp}$.

For $N\geq 5$, we cannot construct the {\em general solution} of equations
(\ref{p2-eqs}) without further assumptions.  However, we can find {\em some}
solutions and these also have this graph symmetry.  Furthermore, if we {\bf
assume} the graph symmetry, then we {\bf can} find the general solution for
some higher values of $N$, but have no general proof that this will be the case
for {\bf all} $N$.

We previously saw this graph symmetry in the context of period $1$ primitives
(see Remark \ref{iota}).
\er   %

\subsection{$5$ Node Quivers of Period 2}

Starting with the general skew-symmetric, $5\times 5$ matrix, with
$$
b_{k1}=m_{k-1},\;\; k=2,\cdots , 5 \quad\mbox{and}\;\; b_{52}=p_1,
$$
we immediately find
$$
b_{32}=m_4 + \varepsilon_{12},\; b_{42}=m_2 + \varepsilon_{14}+\varepsilon(m_1,p_1),\;
b_{43}=m_4 + \varepsilon(m_1,p_1-\varepsilon_{14}),\; b_{53}=p_1 - \varepsilon_{14},\;
b_{54}=m_1,
$$
together with the simple condition $m_3=m_2+\varepsilon_{14}$ and four complicated
conditions.

Imposing the graph symmetry (Remark \ref{graphsymmREM}) leads to
$p_1=m_2=\varepsilon_{14}$, after which two of the four conditions are identically
satisfied, whilst the other pair reduce to a {\em single} condition:
$$
\varepsilon(m_2,p_1)+\varepsilon(m_4,p_1)-\varepsilon_{12}=m_4-m_1 .
$$
We need integer solutions for $m_1,m_2,m_4$.  There are a number of subcases

\subsubsection*{The case $m_1>0,m_4>0$}

In this case the remaining condition reduces to
$$
(|m_2|-m_2-2)(m_1-m_4)=0.
$$
Discarding the period $1$ solution, $m_4=m_1$, we obtain $m_2=-1$, leading to
\bea   %
B(1) &=& \left(\begin{array}{ccccc}
     0 & -m_1 & 1 & 1 & -m_4 \\
     m_1 & 0 & -m_1-m_4 & 1-m_1 & 1 \\
     -1 & m_1+m_4 & 0 & -m_1-m_4 & 1 \\
     -1 & m_1-1  & m_1+m_4 & 0 & -m_1 \\
     m_4 & -1 & -1 & m_1 & 0
     \end{array}\right) , \nn\\
     && \label{n5-1pos4pos}   \\
B(2) &=& \left(\begin{array}{ccccc}
     0 & m_1 & -1 & -1 & m_4 \\
     -m_1 & 0 & -m_4 & 1 & 1 \\
     1 & m_4 & 0 & -m_1-m_4 & 1-m_4 \\
     1 & -1  & m_1+m_4 & 0 & -m_1-m_4 \\
     -m_4 & -1 & m_4-1 & m_1+m_4 & 0
     \end{array}\right) .   \nn
\eea   %
Notice again that $B(2)(m_1,m_4)=\rho B(1)(m_4,m_1) \rho^{-1}$, so the period $2$
property stems from the involution $m_1\leftrightarrow m_4$.

\subsubsection*{The case $m_1>0,\, m_4<0,\, m_2>0$}
There is one condition, which can be reduced by noting that $m_2-m_1m_4>0$,
giving
$$
m_4(m_2-1)=m_1(m_4^2-1).
$$
The left side is {\em negative} and the right {\em positive} unless $m_2=1,\,
m_4=-1$.  We then have $m_3=p_1=m_1+1$, giving
\bea   %
B(1) &=& \left(\begin{array}{ccccc}
     0 & -m_1 & -1 & -m_1-1 & 1 \\
     m_1 & 0 & 1 & -m_1-1 & -m_1-1 \\
     1 & -1 & 0 & 1 & -1 \\
     m_1+1 & m_1+1  & -1 & 0 & -m_1 \\
     -1 & m_1+1 & 1 & m_1 & 0
     \end{array}\right) , \nn\\
     && \label{n5-1pos4neg2pos}   \\
B(2) &=& \left(\begin{array}{ccccc}
     0 & m_1 & 1 & m_1+1 & -1 \\
     -m_1 & 0 & 1 & -m_1-1 & -1 \\
     -1 & -1 & 0 & 1 & 0 \\
     -m_1-1 & m_1+1  & -1 & 0 & 1 \\
      1 & 1 & 0 & -1 & 0
     \end{array}\right) .   \nn
\eea   %

\subsubsection*{The case $m_1>0,\, m_4<0,\, m_2<0$}

Here we have no control over the sign of $m_2-m_1m_4$.

When $m_2-m_1m_4>0$, we have the single condition
$$
(m_2-m_1m_4)(m_2+m_4)+m_1(m_2+1)-m_4 = 0.
$$
Whilst any integer solution would give an example, we have no way of
determining these. (However, Andy Hone has communicated to us that
an algebraic-number-theoretic argument can be used to show that there
are no integer solutions).

When $m_2-m_1m_4<0$, we have $m_4=m_1(m_2+1)$, so
$m_3=m_2-m_1m_4=m_2-m_1^2(m_2+1)$.  We must have $m_2\leq -2$ for $m_4<0$.
Since $m_2-m_1m_4=m_2-m_1^2(m_2+1)<0$, we then choose $m_1$ to be any integer
satisfying $m_1>\sqrt{\frac{m_2}{m_2+1}}$.

Subject to these constraints, the matrices take the form:
\bea   %
B(1) &=& \left(\begin{array}{ccccc}
     0 & -m_1 & -m_2 & -m_3 & -m_1(m_2+1) \\
     m_1 & 0 & -m_1 & m_3(m_1-1) & -m_3 \\
     m_2 & m_1 & 0 & -m_1 & -m_2 \\
     m_3 & -m_3(m_1-1)  & m_1 & 0 & -m_1 \\
     m_1(m_2+1) & m_3 & m_2 & m_1 & 0
     \end{array}\right) , \nn\\
     && \label{n5-1pos4neg2neg}   \\
B(2) &=& \left(\begin{array}{ccccc}
     0 & m_1 & m_2 & m_3 & m_1(m_2+1) \\
     -m_1 & 0 & -m_1(m_2+1) & m_3 & -m_2 \\
     -m_2 & m_1(m_2+1) & 0 & -m_1 & -m_2 \\
     -m_3 & -m_3  & m_1 & 0 & -m_1 \\
      -m_1(m_2+1) & m_2 & m_2 & m_1 & 0
     \end{array}\right) .   \nn
\eea   %
The simplest solution has $m_1=2,\, m_2=-2$.

\subsection{A Family of Period $2$ Solutions}\label{p2-reg}

We are not able to classify all period $2$ quivers. Note that in
Section~\ref{p2} we have classified all sink type period $2$ quivers. In this
section we shall explain how to modify the proof of the classification of
period $1$ quivers (see Section~\ref{p1-gen}) in order to construct a family of
period $2$ quivers (which are, in general, not of sink type).  The introduction
of the involution $\sigma$, defined below, is motivated by the matrices
(\ref{p2n4}) and (\ref{n5-1pos4pos}).

As before, we consider the matrix:
\be  \label{genp2-matrix}    %
B = \left(\begin{array}{ccccc}
      0 & -m_1 &  \cdots & -m_{N-1} \\
      m_1 & 0  && *  \\
      \vdots && 0 &  \\
      m_{N-1} & * && 0
      \end{array}  \right).
\ee    %
However, we assume that, for $r=2,3,\ldots ,N-2$, $m_r=m_{N-r}$ (in the period
$1$ case this property follows automatically). We write $m_{\overline{1}}$
instead of $m_{N-1}$, for convenience. We also assume that $m_1\geq 0$,
$m_{N-1}=m_{\overline{1}}\geq 0$ and $m_1\not=m_{\overline{1}}$ (the last
condition to ensure we obtain strictly period $2$ matrices). We consider the
involution $\sigma$ which fixes $m_r$ for $r\not=1$ and interchanges
$m_1$ and $m_{\overline{1}}$. Let $\mathbf{m}=(m_1,m_2,\ldots
,m_{N-2},m_{\overline{1}})$. We write $\sigma(\mathbf{m})=\sigma(m_1,m_2,\ldots
,m_{N-2}, m_{\overline{1}})=(m_{\overline{1}},m_2,\ldots ,m_{N-2},m_1)$.

Our aim is to construct a matrix $B=B(m_1,m_2,\ldots ,m_{N-1})$
which satisfies the equation
\begin{equation}
\mu_1(B)=\rho B(\sigma(\mathbf{m}))\rho^{-1}.\label{eqn:sigmad}
\end{equation}
Since $\sigma$ is an involution, we shall obtain period $2$ solutions in this
way. As in the period $1$ case, equation~(\ref{eqn:sigmad}) implies that
$(b_{N1},b_{N2},\ldots ,b_{N,N-1})=\sigma(\mathbf{m})$. The derivation of
(\ref{e:mutationequal}) in the period $1$ case is modified by the action of
$\sigma$ to give
\begin{equation} \label{eqn:period2equation}
b_{ij}=\sigma(b_{i-1,j-1})+\varepsilon_{j-1,i-1}.
\end{equation}
An easy induction shows that:
$$
b_{ij}=\sigma^{j-1}(b_{i-j+1,1})+\sum_{s=1}^{j-1}\sigma^{j-1-s}(\varepsilon_{s,i-j+s}).
$$
Applying this in the case $i=N$ we obtain
$$
b_{Nj}=\sigma^{j-1}(b_{N-j+1,1})+\sum_{s=1}^{j-1}\sigma^{j-1-s}(\varepsilon_{s,N-j+s}).
$$
Hence we have
$$
b_{Nj}=\sigma^{j-1}(b_{N-j+1,1})+\varepsilon_{j-1,\overline{1}}+
            \sum_{s=1}^{j-2}\sigma^{j-1-s}(\varepsilon_{s,j-s)})
$$
For $j\leq N-2$, this gives
$$
m_j=m_{N-j}+\sigma^{j-2}(\varepsilon_{1,j-1})+\varepsilon_{j-1,\overline{1}}.
$$
Since $m_j=m_{N-j}$, this is equivalent to
$\sigma^{j-2}(\varepsilon_{1,j-1})+\varepsilon_{j-1,\overline{1}}=0$. For $j=2$ this is
automatically satisfied, since $\varepsilon_{1\overline{1}}=0$.  For $j\geq 3$ and
odd, this is always true.  For $j\geq 4$ and even, this is true if and only if
$m_{j-1}\geq 0$. For $j=N-1$, we obtain
$$
m_1=b_{N,N-1}=\sigma^{N-2}(m_1)+\varepsilon_{2,\overline{1}}+
               \sum_{s=1}^{N-3}\sigma^{N-2-s}(\varepsilon_{s,j-s})
$$
and thus
$$
m_1=\sigma^{N-2}(m_1)+\sigma^{N-3}(\varepsilon_{1,2})+\varepsilon_{2,\overline{1}}.
$$
For $N$ even this is equivalent to
$\varepsilon_{\overline{1},2}+\varepsilon_{2,\overline{1}}=0$, which always holds. For
$N$ odd this gives the condition
$$
m_1=m_{\overline{1}}+\varepsilon_{1,2}+\varepsilon_{2,\overline{1}}.
$$
If $m_2\geq 0$, this is equivalent to $m_1=m_{\overline{1}}$, a contradiction
to our assumption. If $m_2<0$, this is equivalent to
$m_1=m_{\overline{1}}-m_1m_2+m_2m_{\overline{1}}$, which holds if and only if
$m_2=-1$ (since we have assumed that $m_1\not= m_{\overline{1}}$). Therefore,
we obtain a period $2$ solution provided $m_r\geq 0$ for $r$ odd, $r\geq 3$
and, in addition, $m_2=-1$ for $N$ odd.

\section{Recurrences with the Laurent Property}\label{laurent}

As previously said, our original motivation for this work was the well known
connection between cluster algebras and sequences with the Laurent property,
developed by Fomin and Zelevinsky in~\cite{02-3,02-2}.  We note that cluster
algebras were initially introduced (in~\cite{02-3}) in order to study total
positivity of matrices and the (dual of the) canonical basis of
Kashiwara~\cite{91-18} and Lusztig~\cite{90-21} for a quantised enveloping
algebra.

In this section we use the cluster algebras associated to periodic quivers to
construct sequences with the Laurent property. These are likely to be a rich
source of integrable maps.  Indeed, it is well known (see \cite{07-2}) that the
Somos $4$ recurrence can be viewed as an integrable map, having a degenerate
Poisson bracket and first integral, which can be reduced to a $2$ dimensional
symplectic map with first integral.  This $2$ dimensional map is a special case
of the QRT \cite{88-5} family of integrable maps.  The Somos $4$ Poisson
bracket is a special case of that introduced in \cite{03-5} for all cluster
algebra structures.  For many of the maps derived by the construction given in
this section, it is also possible to construct first integrals, often enough to
prove complete integrability.  We do not yet have a complete picture, so do not
discuss this property in general.  However, the maps associated with our
primitives are simple enough to treat in general and can even be linearised.
This is presented in Section \ref{linear}.

A (skew-symmetric, coefficient-free) cluster algebra is an algebraic structure
which can be associated with a quiver. (Recall that we only consider quivers
with no $1$ or $2$-cycles). Given a quiver (with $N$ nodes), we attach a
variable at each node, labelled $(x_1,\cdots ,x_N)$.  When we mutate the quiver
we change the associated matrix according to formula (\ref{gen-mut}) and, {\em
in addition}, we transform the cluster variables $(x_1,\cdots ,x_N)\mapsto
(x_1,\cdots ,\tilde x_\ell,\cdots ,x_N)$, where
\be  \label{ex-rel}   %
x_\ell \tilde x_\ell = \prod_{b_{i\ell}>0} x_i^{b_{i\ell}}+
            \prod_{b_{i\ell}<0} x_i^{-b_{i\ell}},
           \qquad   \tilde x_i = x_i \;\;\mbox{for}\;\; i\neq \ell .
\ee   %
If one of these products is empty (which occurs when all $b_{i\ell}$ have the
same sign) then it is replaced by the number $1$.  This formula is called the
(cluster) {\em exchange relation}.  Notice that it just depends upon the
$\ell^{th}$ column of the matrix.  Since the matrix is skew-symmetric, the
variable $x_\ell$ {\bf does not} occur on the right side of (\ref{ex-rel}).

After this process we have a new quiver $\tilde Q$, with a new matrix $\tilde
B$.  This new quiver has cluster variables $(\tilde x_1,\cdots ,\tilde x_N)$.
However, since the exchange relation (\ref{ex-rel}) acts as the identity on all
except one variable, we write these new cluster variables as $(x_1,\cdots
,\tilde x_\ell,\cdots ,x_N)$.  We can now repeat this process and mutate
$\tilde Q$ at node $p$ and produce a third quiver $\tilde{\tilde Q}$, with
cluster variables $(x_1,\cdots ,\tilde x_\ell,\cdots ,\tilde x_p,\cdots ,x_N)$,
with $\tilde x_p$ being given by an analogous formula (\ref{ex-rel}), but
using variable $\tilde{x}_\ell$ instead of $x_\ell$.

\br[Involutive Property of the Exchange Relation]   %
Since the matrix mutation formula (\ref{gen-mut}) just changes the signs of the
entries in column $n$, a {\bf second} mutation at this node would entail an
{\bf identical} right hand side of (\ref{ex-rel}) (just interchanging the two
products), leading to
$$
\tilde x_\ell \tilde{\tilde x}_\ell = x_\ell \tilde x_\ell
       \quad\Rightarrow\quad     \tilde{\tilde x}_\ell = x_\ell .
$$
Therefore, the exchange relation is an involution.
\er   %

\br[Equivalence of a Quiver and its Opposite]   %
The mutation formula (\ref{ex-rel}) for a quiver and its opposite are identical
since this corresponds to just a change of sign of the matrix entries
$b_{i\ell}$. This is a reason for considering these quivers as {\bf equivalent}
in our context.
\er   %

In this paper we have introduced the notion of {\em mutation periodicity} and
followed the convention that we mutate first at node $1$, then at node $2$,
etc. Mutation periodicity (period $m$) meant that after $m$ steps we return to
a quiver which is equivalent (up to a specific permutation) to the original
quiver $Q$ (see the diagram (\ref{periodchain})).  The significance of this is
that the mutation at node $m+1$ produces an exchange relation which is {\bf
identical in form} (but with a different labelling) to the exchange relation at
node $1$.  The next mutation produces an exchange relation which is {\bf
identical in form} (but with a different labelling) to the exchange relation at
node $2$.  We thus obtain a periodic listing of formulae, which can be
interpreted as an iteration, as can be seen in the examples below.

\subsection{Period $1$ Case}  \label{p1-sequence}

We start with cluster variables $(x_1,\cdots ,x_N)$, with $x_i$ situated at
node $i$.  We then successively mutate at nodes $1,2,3,\dots$ and define
$x_{N+1}=\tilde x_1,\, x_{N+2}=\tilde x_2$, etc.  The exchange relation
(\ref{ex-rel}) gives us a formula of the type
\be   \label{xkxk+n}   %
x_n x_{n+N}=F(x_{n+1},\cdots ,x_{n+N-1}),
\ee  %
with $F$ being the sum of two monomials.  This is interpreted as an $N^{th}$
order recurrence of the real line, with initial conditions $x_i=c_i$ for $i=1,\cdots
,N$. Whilst the right hand side of (\ref{xkxk+n}) is {\em polynomial}, the
formula for $x_{n+N}$ involves a division by $x_n$.  For a general polynomial
$F$, this would mean that $x_n$, for $n>2N$, is a complicated {\bf rational
function} of $c_1,\cdots ,c_N$.  However, in our case, $F$ is derived through
the cluster exchange relation (\ref{ex-rel}), so, by a theorem of \cite{02-3},
$x_n$ is just a {\bf Laurent polynomial} in $c_1,\cdots ,c_N$, for all $n$.  In
particular, if we start with $c_i=1, i=1,\cdots N$, then $x_n$ is an {\bf
integer} for all $n$.

\br[$F$ not $F_n$]  %
For emphasis, we repeat that for a generic quiver we would need to write $F_n$,
since the formula would be different for each mutation.  It is the special
property of period $1$ quivers which enables the formula to be written as a
recurrence.
\er   %

The recurrence corresponding to a general quiver of period $1$ with $N$ nodes
(as described in Theorems~\ref{p1-prototheorem},~\ref{p1-theorem})
corresponding to integers $m_1,m_2,\ldots ,m_{N-1}$ (with $m_r=m_{N-r}$) is:
\be \label{period1recurrence}
x_nx_{n+N}=\prod_{\substack{i=1\\ m_i>0}}^{N-1}x_{n+i}^{m_i}+
\prod_{\substack{i=1\\ m_i<0}}^{N-1}x_{n+i}^{-m_i}
\ee %

\bex[$4$ Node Case]  \label{4node}  {\em   %
Consider Example \ref{n=4p1}, with $m_1=r, m_2=-s$, both $r$ and $s$ positive.
With $r=1, s=2$, the quiver is shown in Figure~\ref{subfig:somos4quiver}. We
start with the matrix
$$
B(1)= \left(\begin{array}{cccc}
     0 & -r & s &  -r \\
     r & 0 & -r(1+s) & s \\
     -s & r(1+s) & 0 & -r \\
     r & -s & r & 0
     \end{array}\right)
$$
and mutate at node $1$, with $(x_1,x_2,x_3,x_4)\mapsto (x_5,x_2,x_3,x_4)$.
Formula (\ref{ex-rel}) gives
\be  \label{x1x5}  %
x_1 x_5 = x_2^r x_4^r+x_3^s ,
\ee  %
whilst the mutation formula (\ref{gen-mut}) gives
$$
B(2) = \left(\begin{array}{cccc}
     0 & r & -s &  r \\
     -r & 0 & -r & s \\
     s & r & 0 & -r(1+s) \\
     -r & -s & r(1+s) & 0
     \end{array}\right).
$$
Note that the second column of this matrix has the same entries (up to
permutation) as the first column of $B(1)$.  This is because $\mu_1 B(1) = \rho
B(1)\rho^{-1}$. Therefore, when we mutate $Q(2)$ at node $2$, with
$(x_5,x_2,x_3,x_4)\mapsto (x_5,x_6,x_3,x_4)$, formula (\ref{ex-rel}) gives
\be   \label{x2x6}  %
x_2 x_6 = x_3^r x_5^r+x_4^s ,
\ee  %
which is of the same form as (\ref{x1x5}), but with indices shifted by $1$.
Formulae (\ref{x1x5}) and (\ref{x2x6}) give us the beginning of the recurrence
(\ref{xkxk+n}), which now explicitly takes the form
$$
x_n x_{n+4} = x_{n+1}^r x_{n+3}^r+x_{n+2}^s .
$$
When $r=1, s=2$, this is exactly the Somos $4$ sequence (\ref{somos4}).  When
$r=s=1$, we obtain the recurrence considered by Dana Scott (see \cite{91-17}
and \cite{08-3}).  This case was also considered by Hone (see Theorem 1 in
\cite{07-2}), who showed that it is {\em super-integrable} and linearisable.
}\eex   %

\bex[$5$ Node Case]  \label{5node1}  {\em   %
Consider Example \ref{n=5p1}, with $m_1=r, m_2=-s$, both $r$ and $s$ positive.
We start with the matrix
$$
B = \left(\begin{array}{ccccc}
     0 & -r & s & s & -r \\
     r & 0 & -r(1+s) & -s(r-1) & s \\
     -s & r(1+s)  & 0 & -r(1+s)  & s \\
     -s & s(r-1)  & r(1+s) & 0 & -r \\
     r & -s & -s & r & 0
     \end{array}\right)
$$
and mutate at node $1$, with $(x_1,x_2,x_3,x_4,x_5)\mapsto
(x_6,x_2,x_3,x_4,x_5)$. Formula (\ref{ex-rel}) gives
$$
x_1 x_6 = x_2^r x_5^r+x_3^s x_4^s .
$$
Proceeding as before, the general term in the recurrence (\ref{xkxk+n}) takes
the form
$$
x_n x_{n+5} = x_{n+1}^r x_{n+4}^r+x_{n+2}^s x_{n+3}^s,
$$
which reduces to Somos $5$ when $r=s=1$ (giving us the quiver of Figure
\ref{fig:somos5quiver}).
}\eex   %

\bex[$6$ Node Case]  \label{6node1}  {\em   %
Consider Example \ref{n=6p1}.  The first thing to note is that there are $3$
parameters $m_i$, so we have rather more possibilities in our choice of signs.
Having already obtained Somos $4$ and Somos $5$, one may be lured into thinking
that Somos $6$ will arise.  However, Somos $6$
$$
x_n x_{n+6} = x_{n+1} x_{n+5}+x_{n+2} x_{n+4}+x_{n+3}^2
$$
has $3$ terms, so cannot directly arise through the cluster exchange relation
(\ref{ex-rel}), although we remark that it is shown in \cite{02-2} that the
terms in the Somos $6$ and Somos $7$ sequences are Laurent polynomials in their
initial terms.  However, various subcases of Somos $6$ {\bf do} arise in our
construction.  They are, in fact, special cases of the Gale-Robinson sequence
of Example \ref{galerob}.

\paragraph{The case $m_1=r,m_2=-s,m_3=0$ with $r, s$ positive.}  We can read
off the recurrence from the first column of the matrix of Example \ref{n=6p1}, which
is $(0,r,-s,0,-s,r)^T$, giving
$$
x_n x_{n+6} = x_{n+1}^r x_{n+5}^r+x_{n+2}^s x_{n+4}^s,
$$
which gives the first two terms of Somos $6$ when $r=s=1$.

\paragraph{The case $m_1=r,m_2=0,m_3=-s$ with $r, s$ positive.}  The first
column of the matrix is now $(0,r,0,-s,0,r)^T$, giving
$$
x_n x_{n+6} = x_{n+1}^r x_{n+5}^r+x_{n+3}^s .
$$
For a subcase of Somos $6$ we choose $r=1, s=2$.

\paragraph{The case $m_1=0,m_2=r,m_3=-s$ with $r, s$ positive.}  The first
column of the matrix is now $(0,0,r,-s,r,0)^T$, giving
$$
x_n x_{n+6} = x_{n+2}^r x_{n+4}^r+x_{n+3}^s ,
$$
again with $r=1, s=2$.
}\eex   %

\bex[Gale-Robinson Sequence ($N$ nodes)]\label{galerob}  {\em   %
The $2-$term Gale-Robinson recurrence (see Equation (6) of \cite{91-17}) is
given by
$$
x_nx_{n+N}=x_{n+N-r}x_{n+r}+x_{n+N-s}x_{n+s} ,
$$
for $0<r<s\leq N/2$, and is one of the examples highlighted in \cite{02-2}.  We
remark that this corresponds to the period $1$ quiver with $m_r=1$ and $m_s=-1$
(unless $N=2s$, in which case we take $m_s=-2$); see Theorem~\ref{p1-theorem}.
}\eex  %

\subsection{Period $2$ Case}
\label{ss:period2case}

We start with cluster variables $(z_1,\cdots ,z_N)$, with $z_i$ situated at
node $i$.  We then successively mutate at nodes $1,2,3,\dots$ and define
$z_{N+1}=\tilde z_1,\, z_{N+2}=\tilde z_2$, etc.  However, the exchange
relation (\ref{ex-rel}) now gives us an {\em alternating pair} of formulae of
the type
\be  \label{f0f1}   %
z_{2n-1} z_{2n-1+N} = F_0(z_{2n},\cdots ,z_{2n+N-2}),\;\;
    z_{2n} z_{2n+N} = F_1(z_{2n+1},\cdots ,z_{2n+N-1}), \;\; n=1,2,\cdots.
\ee  %
with $F_i$ being the sum of two monomials.  It is natural, therefore, to
relabel the cluster variables as $x_n=z_{2n-1},\, y_n=z_{2n}$ and to
interpret~\eqref{f0f1} as a two-dimensional recurrence for $(x_n,y_n)$.
When $N=2m$, the recurrence is of order $m$.
When $N=2m-1$, the recurrence is again of order $m$, but the first exchange relation
plays the role of a {\bf boundary condition}.  We need $m$ points in the plane
to act as initial conditions.  When $N=2m$, the values $z_1,\cdots ,z_{2m}$
define these $m$ points.  When $N=2m-1$, we need $z_{2m}$ {\bf in addition} to
the {\bf given} initial conditions $z_1,\cdots ,z_{2m-1}$.  Again, since our
recurrences are derived through the cluster exchange relation (\ref{ex-rel}), the
formulae for $(x_n,y_n)$ are Laurent polynomials of initial conditions.  In the
case of $N=2m-1$, this really does mean initial conditions $z_1,\cdots
,z_{2m-1}$.  The expression for $y_m=z_{2m}$ is already a polynomial, so it is
important that it {\bf does not} occur in the denominators of later terms.

\bex[$4$ Node Case]  \label{4nodep2}  {\em   %
Consider the general period $2$ quiver with $4$ nodes, which has corresponding
matrices (\ref{p2n4}), which we write with $m_1=r, m_2=-s, m_3=t$, where
$r,s,t$ are positive:
\bea     %
B(1) &=& \left(\begin{array}{cccc}
     0 & -r & s & -t \\
     r & 0 & -t-rs & s \\
     -s & t+rs & 0 & -r \\
     t & -s & r & 0
     \end{array}\right) ,   \nn\\
     && \label{p2n4rst}   \\  %
B(2) &=& \left(\begin{array}{cccc}
     0 & r & -s & t \\
     -r & 0 & -t & s \\
     s & t & 0 & -r-st \\
     -t & -s & r+st & 0
     \end{array}\right) .  \nn
\eea    %
Mutating $Q(1)$ at node $1$, with $(z_1,z_2,z_3,z_4)\mapsto (z_5,z_2,z_3,z_4)$,
formula (\ref{ex-rel}) gives
\be  \label{z1z5}  %
z_1 z_5 = z_2^r z_4^t+z_3^s ,
\ee  %
whilst mutating $Q(2)$ at node $2$, with $(z_5,z_2,z_3,z_4)\mapsto
(z_5,z_6,z_3,z_4)$, formula (\ref{ex-rel}) gives
\be  \label{z2z6}  %
z_2 z_6 = z_3^t z_5^r+z_4^s ,
\ee  %
When $t\neq r$ these formulae are not related by a shift of index.  However,
since $B(3)=\mu_2 B(2) = \rho^2 B(1)\rho^{-2}$, mutating $Q(3)$ at node $3$, with
$(z_5,z_6,z_3,z_4)\mapsto (z_5,z_6,z_7,z_4)$, leads to
\be  \label{z3z7}  %
z_3 z_7 = z_4^r z_6^t+z_5^s ,
\ee  %
which is just (\ref{z1z5}) with a shift of $2$ on the indices.  This pattern
continues, giving
\be  \label{xnyn}  %
x_n x_{n+2} = y_n^r y_{n+1}^t+x_{n+1}^s ,\quad
        y_n y_{n+2} = x_{n+1}^t x_{n+2}^r + y_{n+1}^s .
\ee   %
The appearance of $x_{n+2}$ in the definition of $y_{n+2}$ is not a problem,
since it can be replaced by the expression given by the first equation.

As shown in Figure \ref{p2quivers}, we could equally start with the matrices
$$
\bar B(1) = \rho^{-1} B(2) \rho , \quad  \bar B(2) = \rho B(1) \rho^{-1}.
$$
Since $\bar B(1)(r,s,t)=B(1)(t,s,r), \bar B(2)(r,s,t)=B(2)(t,s,r)$, we obtain a
two-dimensional recurrence
\be  \label{unvn}  %
u_n u_{n+2} = v_n^t v_{n+1}^r+u_{n+1}^s ,\quad
        v_n v_{n+2} = u_{n+1}^r u_{n+2}^t + v_{n+1}^s ,
\ee   %
where we have labelled the nodes as $\zeta_1, \zeta_2,\cdots$ and then
substituted $u_k=\zeta_{2k-1}, v_k = \zeta_{2k}$.  With initial conditions
$(z_1,z_2,z_3,z_4)=(1,1,1,1)$ and
$(\zeta_1,\zeta_2,\zeta_3,\zeta_4)=(1,1,1,1)$, recurrences (\ref{xnyn}) and
(\ref{unvn}) generate {\em different} sequences of integers.  However, just
making the change $\zeta_4=2$, reproduces the original $z_n$ sequence.  This
corresponds to a shift in the labelling of the nodes, given by
$$
u_n = y_n, \quad v_n= x_{n+1}, \;\;\; n=1,2,\cdots
$$
}\eex   %

\bex[$5$ Node Case]  \label{5nodep2}  {\em   %
Consider the case with matrices (\ref{n5-1pos4pos}), which we write with
$m_1=r, m_4=t$, where $r,t$ are positive.  The same procedure leads to the recurrence
\be  \label{xnyn5node}  %
 y_n x_{n+3} = y_{n+2}^r x_{n+1}^t+x_{n+2} y_{n+1},  \quad
   x_{n+1} y_{n+3} = y_{n+1}^r x_{n+3}^t + x_{n+2} y_{n+2}, \;\;\; n=1,2,\cdots
\ee   %
together with
$$
x_1 y_3=y_1^r x_3^t + x_2 y_2,
$$
and initial conditions $(x_1,y_1,x_2,y_2,x_3)=(c_1,c_2,c_3,c_4,c_5)$.  The
iteration (\ref{xnyn5node}) is a third order two-dimensional recurrence and $y_3$ acts as
the sixth initial condition.

As above, it is possible to construct a companion recurrence, corresponding to the
choice
$$
\bar B(1) = \rho^{-1} B(2) \rho , \quad  \bar B(2) = \rho B(1) \rho^{-1}.
$$
}\eex   %

\section{Linearisable Recurrences From Primitives}\label{linear}

This section is concerned with the recurrences derived from {\em period $1$
primitives}. Similar results can be shown for higher periods, but we omit these
here.

Our primitive quivers are inherently simpler than composite ones (as their name
suggests!). The mutation process (at node $1$) reduces to a simple matrix
conjugation.  The cluster exchange relation is still nonlinear, but turns out
to be {\bf linearisable}, as is shown in this section.

Consider the $k^{th}$ (period $1$) primitive $P_N^{(k)}$ with $N$ nodes, such
as those depicted in Figures~\ref{234node} to \ref{6node}.  As before, we
attach a variable at each node, labelled $(x_1,\cdots ,x_N)$, with $x_i$
situated at node $i$ for each $i$.
We then successively mutate at nodes $1,2,3,\dots$ and
define $x_{N+1}=\tilde x_1,\, x_{N+2}=\tilde x_2$, etc.
At the $n$th mutation, we start with the cluster
$\{x_n,x_{n+1}, \ldots ,x_{N+n-1}\}$.
By the periodicity property, the corresponding quiver is always $P_N^{(k)}$.
The exchange relation (\ref{ex-rel}) gives us the formula
\be   \label{en}   %
x_n x_{n+N}=x_{n+k}x_{n+N-k}+1,
\ee  %
where $x_{n+N}$ is the new cluster variable replacing $x_n$.
Note that one of the products in (\ref{ex-rel}) is empty. This is the $n^{th}$
iteration, which we label $E_n$. For $gcd(k,N)=1$, this is a genuinely new sequence
for each $N$. However, when $gcd(k,N)=m>1$, the sequence (\ref{en})
decouples into $m$ copies of an iteration of order $(N/m)$.

Specifically, if $N=ms$ and $k=mt$, for integers $s,t$, the quiver
$P_N^{(k)}$ separates into $m$ disconnected components (see Figures \ref{subfig:P42},
\ref{subfig:P62} and \ref{subfig:P63}). The corresponding sequence decouples
into $m$ copies of the sequence associated with the primitive $P_s^{(t)}$,
since (\ref{en}) then gives
$$
x_nx_{n+ms}=x_{n+mt}x_{n+(s-t)m}+1 .
$$
With $n=ml+r,\, y_l^{(r)}=x_{ml+r},\;\; 0\leq r\leq m-1$, this gives $m$ identical
iteration formulae
\be  \label{k-1-reduce}  %
y_l^{(r)}y_{l+s}^{(r)}=y_{l+t}^{(r)}y_{l+s-t}^{(r)}+1 .
\ee  %
Thus if, in (\ref{en}), we use the initial conditions $x_i=1,\, 1\leq i\leq N$,
we obtain $m$ copies of the integer sequence generated by (\ref{k-1-reduce}).

\subsection{First Integrals}

Subtracting the two equations $E_n$ and $E_{n+k}$ (see~\eqref{en})
leads to
$$
\frac{x_n+x_{n+2k}}{x_{n+k}} = \frac{x_{n+N-k}+x_{n+N+k}}{x_{n+N}}.
$$
With the definition
\be  \label{jnk}   %
J_{n,k}=\frac{x_n+x_{n+2k}}{x_{n+k}},
\ee  %
we therefore have
\be  \label{periodj}  %
J_{n+N-k,k}=J_{n,k},
\ee  %
giving us $N-k$ independent functions $\{J_{i,k}:1\leq i\leq N-k\}$ (or
equivalently $\{J_{i,k}:n\leq i\leq n+N-k-1\}$).

\br[Decoupled case]  %
Again, when $gcd(N,k)=m>1$, the sequence (\ref{en}) decouples
into $m$ copies of  (\ref{k-1-reduce}) and the sequence $J_{n,k}$ (with
periodicity $N-k$) splits into $m$ copies of the corresponding sequence of
$J$'s for the primitive $P_s^{(t)}$ (where $N=ms$ and $k=mt$), since, putting
$n=ml+r$ and $I_{l,t}^{(r)}=J_{ml+r,k}$ we obtain
$$I_{l,t}^{(r)}=\frac{x_{ml+r}+x_{ml+r+2mt}}{x_{ml+r+mt}}=
\frac{y_l^{(r)}+y_{l+2t}^{(r)}}{y_{l+t}^{(r)}},$$
satisfying $I_{l+s-t,t}^{(r)}=I_{l,t}^{(r)}$.
\er   %

Let $\alpha$ be any function of $N-k$ variables and define
$\alpha^{(n)}=\alpha(J_{n,k},\cdots ,J_{n+N-k-1,k})$.  Then, from the
periodicity (\ref{periodj}), $\alpha^{(n+N-k)}=\alpha^{(n)}$ (it can happen
that the function will have periodicity $r\leq N-k$). Then the function
$$
K_\alpha^{(n)} = \sum_{i=0}^{N-k-1} \alpha^{(n+i)}
$$
is a first integral for the recurrence (\ref{en}), meaning that it satisfies
$K_\alpha^{(n+1)}=K_\alpha^{(n)}$.  It is thus always possible to construct,
for the recurrence (\ref{en}), $N-k$ independent first integrals.  For $k=1$ this is
the maximal number of integrals, unless the recurrence is itself periodic (see
\cite{91-4} for the general theory of integrable maps).

For example, $N-k$ independent first integrals $\{K_p^{(n)}:1\leq p\leq N-k\}$
are given by
$$
    K_p^{(n)} = \sum_{i=0}^{N-k-1} \alpha_p^{(n+i)}, \quad\mbox{where}\quad
      \alpha_p^{(n)} = \prod_{i=0}^{p-1} J_{n+i,k}.
$$
Using the condition (\ref{periodj}) and the definition (\ref{jnk}), it can be
seen that the $\alpha_p^{(n)}$ depend upon the variables $x_n,\cdots ,x_{n+N+k-1}$, so
equation (\ref{en}) must be used to eliminate $x_{n+N},\cdots ,x_{n+N+k-1}$ in
order to get the correct form of these integrals in terms of the $N$
independent coordinates.

\bex
As an example, consider the case $N=4$ and $k=1$. This corresponds to
the recurrence:
\be
x_nx_{n+4}=x_{n+1}x_{n+3}+1
\label{e:egrecurrence}
\ee
for the primitive $P_4^{(1)}$. We have
$$J_{n,1}=\frac{x_n+x_{n+2}}{x_{n+1}}.$$
Then $\alpha_1^{(n)}=J_{n,1}$, $\alpha_2^{(n)}=J_{n,1}J_{n+1,1}$ and
$\alpha_3^{(n)}=J_{n,1}J_{n+1,1}J_{n+2,1}$.
So
\begin{align*}
K_1^{(n)} &= \alpha_1^{(n)}+\alpha_1^{(n+1)}+\alpha_1^{(n+2)}=
J_{n,1}+J_{n+1,1}+J_{n+2,1}; \\
K_2^{(n)} &= \alpha_2^{(n)}+\alpha_2^{(n+1)}+\alpha_2^{(n+2)}=
J_{n,1}J_{n+1,1}+J_{n+1,1}J_{n+2,1}+J_{n+2,1}J_{n+3,1} \\
&= J_{n,1}J_{n+1,1}+J_{n+1,1}J_{n+2,1}+J_{n+2,1}J_{n,1}; \\
K_3^{(n)} &= \alpha_3^{(n)}+\alpha_3^{(n+1)}+\alpha_3^{(n+2)}
=J_{n,1}J_{n+1,1}J_{n+2,1}+J_{n+1,1}J_{n+2,1}J_{n+3,1}+J_{n+2,1}J_{n+3,1}J_{n+4,1} \\
&=3J_{n,1}J_{n+1,1}J_{n+2,1}.
\end{align*}
Using~(\ref{e:egrecurrence}), we obtain $N-k=3$ independent first integrals. For
simplicity we write $a=x_n$, $b=x_{n+1}$, $c=x_{n+2}$ and $d=x_{n+3}$:
\begin{align*}
K_1^{(n)} &= \frac{a}{b}+\frac{b}{a}+\frac{b}{c}+\frac{c}{b}+\frac{c}{d}+\frac{d}{c}+\frac{1}{ad}; \\
K_2^{(n)} &= 3+\frac{a}{c}+\frac{c}{a}+\frac{b}{d}+\frac{d}{b}+\frac{1}{ac}+\frac{1}{bd}+
\frac{ad}{bc}+\frac{bd}{ac}+\frac{ac}{bd}+\frac{b^2}{ac}+\frac{c^2}{bd}+\frac{b}{acd}+\frac{c}{abd}; \\
K_3^{(n)} &= 3(\frac{a}{b}+\frac{b}{a}+\frac{b}{c}+\frac{c}{b}+\frac{c}{d}+\frac{d}{c}+
\frac{a}{d}+\frac{d}{a}+\frac{1}{ab}+\frac{1}{bc}+\frac{1}{cd}+\frac{1}{ad}).
\end{align*}
\eex

\br[Decoupled case]  %
Again, when $gcd(N,k)>1$, the sequence (\ref{en}) decouples into
$m$ copies of  (\ref{k-1-reduce}) and we use the first integrals built out of
the functions $I_{l,t}^{(r)}$.
\er   %

Let the sequence $\{x_n\}$ be given by the iteration (\ref{en}), with initial
conditions $\{x_i=a_i: 1\leq i\leq N\}$.  We have $K_p^{(n)}=K_p^{(1)}$, which
is evaluated in terms of $a_i$.  We also have $\{J_{i,k}=c_i:1\leq i\leq
N-k\}$, together with the periodicity condition (\ref{periodj}), which can also
be written as $J_{n,k}=c_n$ with $c_{n+N-k}=c_n$.  The first integrals
$K_p^{(n)}$ have simpler formulae when written in terms of $c_1,\cdots
,c_{N-k}$ (each of which is a rational function of the $a_i$).

\br[Complete integrability] The complete integrability of the maps associated with the
$P_N^{(1)}$ ($N$ even) is shown in~\cite{Fordypreprint}.
\er

\subsection{A Linear Difference Equation}\label{lin:diff:eq}

We show in this subsection that the difference equation (\ref{en}) can be
linearised.

\bt[Linearisation]  \label{linear-thm}  %
If the sequence $\{x_n\}$ is given by the iteration (\ref{en}), with initial
conditions $\{x_i=a_i: 1\leq i\leq N\}$, then it also satisfies
\be  \label{lineq}  %
x_n+x_{n+2k(N-k)}=S_{N,k} x_{n+k(N-k)},
\ee  %
where $S_{N,k}$ is a function of $c_1,\cdots ,c_{N-k}$, which is symmetric
under cyclic permutations.
\et  %

\noindent \textbf{Proof of case $k=1$}: We first prove this theorem for the case
$k=1$, later showing that the general case can be reduced to this.

We fix $k=1$. For $i\in \mathbb{N}$, let $L_i=x_i+x_{i+2}-c_ix_{i+1}$. For
$1\leq i\leq 2N-3$, we have that $J_{i,1}=c_i$ (see the last paragraph of the
previous section), from which it follows that $L_i=0$, but we regard the $x_i$
as formal variables for the time being (see the end of the proof of
Proposition~\ref{prop:linearisation}). For $i=0,1,\ldots ,2N-2$, we define a
sequence $a_i$ as follows. Set $a_0=0, a_1=1$ and then, for $2\leq n\leq N-1$,
define $a_n$ recursively by:
\begin{equation} \label{e:anrecursion}
a_n=-a_{n-2}-c_{n-1}a_{n-1}.
\end{equation}
We also set $b_{2N-2}=0$, $b_{2N-3}=1$ and then, for $N-1\leq n\leq 2N-3$,
define $b_n$ recursively by:
\begin{equation} \label{e:bnrecursion}
b_n=-b_{n+2}-c_{n+1}b_{n+1}.
\end{equation}

\bl \label{lem:abflip} For $0\leq n\leq N-1$, we have
$b_{2N-2-n}=a_n|_{c_l\mapsto c_{2N-2-l}}$. \el

\noindent \textbf{Proof}: This is easily shown using induction on $n$ and
equations~(\ref{e:anrecursion}) and~(\ref{e:bnrecursion}).$\Box$

The proofs of the following results (Lemma~\ref{lem:andescription},
Proposition~\ref{prop:linearisation} and
Corollary~\ref{cor:generallinearisation}) will be given in the Appendix. We
first describe the $a_n$ explicitly. Define:
\begin{align*}
t_{k,\text{odd}}^n &= \sum_{ \substack{
1\leq i_1<i_2<\cdots <i_k\leq n \\
{i_1\text{\ odd},i_2\text{\ even},\ldots }}}
c_{i_1}c_{i_2}\cdots c_{i_k} \\
t_{k,\text{even}}^n &= \sum_{ \substack{
1\leq i_1<i_2<\cdots <i_k\leq n \\
{i_1\text{\ even},i_2\text{\ odd},\ldots }}} c_{i_1}c_{i_2}\cdots c_{i_k}
\end{align*}

\bl \label{lem:andescription}
Suppose that $0\leq n\leq N-1$. Then \\
(a) If $n=2r$ is even,
$$
a_{2r}=(-1)^r \sum_{k=0}^{r-1} (-1)^k t_{2k+1,\text{odd}}^{(2r-1)}.
$$
(b) If $n=2r-1$ is odd,
$$
a_{2r-1}=(-1)^{r-1} \sum_{k=0}^{r-1} (-1)^k t_{2k,\text{odd}}^{(2r-2)}.
$$
(c) We have $a_n=a_n|_{c_l\mapsto c_{n-l}}$ and $a_{N-1}=b_{N-1}$.
\el  %

For $n\in\mathbb{N}$ and $0\leq k\leq n$, define:
$$t_{k,alt}^n =
\begin{cases}
\sum_{ \substack{
1\leq i_1<i_2<\cdots <i_k\leq n \\
i_1,i_2,\ldots ,i_k\text{ of alternating parity}}} c_{i_1}c_{i_2}\cdots c_{i_k}
& \text{ if }k>0; \\
2 & \text{ if }k=0.
\end{cases}$$

Let $L=\sum_{i=1}^{N-1} (-1)^i a_i L_i + \sum_{i=N}^{2N-3} (-1)^i b_i L_i$.
Since $a_1=b_{2N-3}=1$, the coefficients of $x_1$ and $x_{2N-1}$ in $L$ are
both $1$. By equations~(\ref{e:anrecursion}) and~(\ref{e:bnrecursion}), the
coefficient of $x_i$ in $L$ is zero for $i=2,3,\ldots ,N-1,N+2,\ldots ,2N-2$.
By Lemma~\ref{lem:andescription}(c), $a_{N-1}=b_{N-1}$, and it follows that the
coefficient of $x_{N+1}$ is also zero. The coefficient of $x_N$ is
$$S_{N,1}=(-1)^{N-2}(a_{N-2}+c_{N-1}a_{N-1}+b_N).$$
Note that $b_N=a_{N-2}|_{c_l\mapsto c_{2N-2-l}}$ by Lemma~\ref{lem:abflip}, so
$b_N=a_{N-2}|_{c_l\mapsto c_{N-1-l}}$, since $c_{n+N-1}=c_n$. This allows us to
compute the coefficient of $x_N$ explicitly to give:

\bp \label{prop:linearisation} We have
$$
x_1+x_{2N-1}=S_{N,1}x_N,
$$
where
$$S_{N,1}=
\begin{cases}
(-1)^{r-1} \sum_{k=0}^{r-1} (-1)^k t_{2k+1,\text{alt}}^{(2r-1)} &
\text{\ if\ }N=2r\text{\ is\ even}; \\
(-1)^{r-1} \sum_{k=0}^{r-1} (-1)^k t_{2k,\text{alt}}^{(2r-2)} & \text{\ if\
}N=2r-1\text{\ is\ odd}.
\end{cases}
$$
\ep

\bc \label{cor:generallinearisation} For all $n\in\mathbb{N}$,
$$
x_n+x_{n+2(N-1)}=S_{N,1}x_{n+N-1},
$$
where $S_{N,1}$ is as above.
\ec  %

\bex   %
We calculate $S_{N,1}$ for some small values of $N$. We have:
\begin{align*}
S_{2,1} &= c_1; \\
S_{3,1} &= c_1c_2-2; \\
S_{4,1} &= c_1c_2c_3-c_1-c_2-c_3; \\
S_{5,1} &= c_1c_2c_3c_4-c_1c_2-c_2c_3-c_3c_4-c_4c_1+2; \\
S_{6,1} &= c_1c_2c_3c_4c_5-c_1c_2c_3-c_2c_3c_4-c_3c_4c_5-c_4c_5c_1-c_5c_1c_2+c_1+c_2+c_3+c_4+c_5; \\
S_{7,1} &= c_1c_2c_3c_4c_5c_6-c_1c_2c_3c_4-c_2c_3c_4c_5-c_3c_4c_5c_6-c_4c_5c_6c_1-c_5c_6c_1c_2-c_6c_1c_2c_3 \\
&\quad\quad +c_1c_2+c_1c_4+c_1c_6+c_3c_4+c_3c_6+c_5c_6+c_2c_3+c_2c_5+c_4c_5-2.
\end{align*}
We remark that $N=7$ gives the first example where the terms of fixed degree in
$S_{N,1}$ (in this case degree $2$) do not form a single orbit under the cyclic
permutation $(1\ 2\ \cdots \ N-1)$.
\eex  %

The case of $N=4$ can be found in \cite{07-2}.

\subsubsection{The case of general $k>1$}

When $k>1$, the system of equations $J_{n,k}=c_n$ (with $c_{n+N-k}=c_n$) splits
into $k$ subsystems. Writing $n=mk-r$ for some $m\geq 1$ and $0\leq r<k$ we
define
$$
z_m=x_{mk-r}, \quad\mbox{and}\quad I_{m,1}^{(r)} = \frac{z_m+z_{m+2}}{z_{m+1}}.
$$
Writing $J_{n,k}$ (see (\ref{jnk})) in terms of $z_m$, we see that
$J_{n,k}=I_{m,1}^{(r)}$. Define $M=N-k+1$, so $c_{n+M-1}=c_n$.  If
$gcd(N,k)=1$, then, for each $r$, $I_{m,1}^{(r)}$ cycle through {\bf all} of
$c_1,\dots ,c_{M-1}$ (in some order).  For $r=k-1$, label this sequence of
$c_i$ as $d_1,\dots ,d_{M-1}$.  It is important to note that, for other values
of $r$, the order is just a cyclic permutation of $d_1,\dots ,d_{M-1}$.  We
therefore have the conditions for Corollary \ref{cor:generallinearisation},
giving
$$
z_m+z_{m+2(M-1)}=S_{M,1}(d_1,\dots ,d_{M-1}) z_{m+(M-1)}.
$$
Writing this in terms of $x_n$ gives (\ref{lineq}) with
$S_{N,k}=S_{M,1}(d_1,\dots ,d_{M-1})$, given by Proposition
\ref{prop:linearisation}.

When $(N,k)\neq 1$, we should first use (\ref{k-1-reduce}) to reduce to the
relatively prime case and proceed as above.

\br   %
We need $2k(N-k)$ initial conditions in order to generate a sequence with
(\ref{lineq}), but are only supplied with $\{x_i=a_i: 1\leq i\leq N\}$.  If we
use the iteration (\ref{en}) to generate the remaining initial conditions for
(\ref{lineq}), then (\ref{en}) and (\ref{lineq}) will generate exactly the same
sequence of numbers.
\er  %

\subsection{Pell's Equation}

For $k=1$, the sequence (\ref{en}) arising from the primitive $B_n^1$ has
entries which are closely related to Pell's equation, as indicated to us by
examples in~\cite{09-1}, e.g. sequences A001519 and A001075 for $N=2,N=3$
respectively. By Theorem \ref{linear-thm}, we have
\begin{equation} \label{e:linearprimitive1}
x_n+x_{n+2(N-1)}=S_{N,1}x_{n+N-1},
\end{equation}
for $n\geq 1$. We have set $x_n=1$ for
$1\leq n\leq N$, and it is easy to check that $x_n=n-N+1$ for
$N\leq n\leq 2N-1$. It follows that $S_{N,1}=N+1$.
Subsequences of the form $y_m=x_{m(N-1)+c}$ for some constant $c$ satisfy
the recurrence $y_m+y_{m+2}=(N+1)y_{m+1}$ which has associated quadratic
equation $\lambda^2-(N+1)\lambda+1=0$, with roots
\begin{equation} \label{e:pellroots}
\alpha_{\pm}=\frac{N+1\pm\sqrt{(N+1)^2-4}}{2}.
\end{equation}

\begin{prop} \ \\
(a) Suppose that $N=2r-1$ is odd. For $m\in\mathbb{Z}$, $m\geq 0$, let
$a_m=x_{(N-1)m+r}$. Choose $1\leq t\leq N-1$, and let
$b_m=x_{(N-1)m+t+1}-x_{(N-1)m+t}$. Then the pairs $(a_m,b_m)$ for $m>0$ are
the positive integer solutions of the Pell equation $a^2-(r^2-1)b^2=1$. \\
(b) Suppose that $N=2r$ is even. Choose $t,t'$ such that $1\leq t\leq r$ and
$1\leq t'\leq N-1$. For $m\in\mathbb{Z}$, $m\geq 0$, let
$a_m=x_{(N-1)m+t}+x_{(N-1)m+N+1-t}$ and let
$b_m=x_{(N-1)m+t'+1}-x_{(N-1)m+t'}$. Then the pairs $(a_m,b_m)$ for $m>0$ are
the positive integer solutions of the Pell equation $a^2-((2r+1)^2-4)b^2=4$.
\end{prop}

\textbf{Proof}:
The general solution of $y_m+y_{m+2}=(N+1)y_{m+1}$ is
$$
y_m=A_+\alpha_+^{m-1}+A_-\alpha_-^{m-1}
$$
for arbitrary constants $A_{\pm}$. The description above of the initial
terms in the sequence $(x_n)$ gives initial terms (for $m=0$ and $1$) for
the subsequences $a_m$ and $b_m$ in each case, and it follows that, in the
odd case, $$a_m+b_m\sqrt{r^2-1}=(r+\sqrt{r^2-1})^m,$$
and, in the even case,
$$
a_m+b_m\sqrt{(2r+1)^2-4}=2^{1-m}(2r+1+\sqrt{(2r+1)^2-4})^m.
$$
In the odd case, it is well-known that these are the positive integer solutions
to $a^2-(r^2-1)b^2=1$, and in the even case, the description of the solutions
is given in~\cite{78-11} (see also~\cite[Theorem 1]{98-7}). (For the $N=2$
case, see for example~\cite{09-1}, sequence A001519). $\Box$

\section{Parameters and Coefficients}\label{s:coefficients}

We recalled the definition of a skew-symmetric coefficient-free cluster algebra
in Section~\ref{laurent}. The general definition~\cite{02-3} of a cluster
algebra allows for coefficients in the exchange relations. We use the \emph{ice
quiver} approach of~\cite[2.2]{07-5} in which some of the cluster variables are
specified to be frozen. The definition of the cluster algebra is the same,
except that mutation at the frozen cluster variables is not allowed.
Our aim in this section is to describe the period $1$ ice quivers.
Each such quiver models a corresponding Laurent recurrence with parameters,
again via the Laurent phenomenon~\cite[3.1]{02-3}.
In other words, we will give an answer to the question as to when can we
take an iterative binomial recurrence coming from a periodic quiver
and add coefficients to the recurrence and still explain this
recurrence in terms of periodic mutations of a frozen quiver using our
methods. We will then give some examples of recurrences modelled by periodic
ice quivers.

We consider an initial cluster consisting of $N$ unfrozen cluster variables
$x_1,x_2,\ldots ,x_N$ and $M$ frozen cluster variables $y_1,y_2,\ldots ,y_M$.
Thus, each seed contains a cluster with $N$ unfrozen cluster variables together
with the frozen variables $y_1,\ldots ,y_M$, which never change. The quiver in
the seed has $N$ unfrozen vertices $1,2,\ldots ,N$ and $M$ frozen vertices
$N+1,\ldots ,N+M$. The exchange matrix $B$ will be taken to be the
corresponding skew-symmetric matrix. The entries $b_{N+i,N+j},\, 1\leq i,j\leq
M$ do not play a role, so we take them to be zero (equivalently, there are no
arrows from a vertex $N+1,\ldots ,N+M$ of the quiver to another such vertex).

Note that in the usual frozen variable set-up, columns $N+1\ldots , N+M$ of $B$
are not included. This makes no difference, since the entries in these columns
do not appear in the exchange relations. They are determined by the rest of $B$
since $B$ is skew-symmetric, and by the above assumption on zero entries. In
order to ensure that the entries $b_{N+i,N+j}$ remain zero, we must modify the mutation
$\mu_i$ slightly: $\widetilde{\mu_i}$ is the same as $\mu_i$ except that the
entries $b_{N+i,N+j}$, $1\leq i,j\leq M$, remain zero by definition. We find it
convenient to include the extra columns in order to study the period $1$ ice
quiver case.

The exchange relation can then be written as follows, for $1\leq \ell\leq n$:

\be  \label{e:coefficientexchangerelation}  %
x_\ell \tilde x_\ell = \mathop{\prod_{1\leq i\leq M}}_{b_{N+i,\ell}>0}
y_i^{b_{N+i,\ell}} \mathop{\prod_{1\leq i\leq N}}_{b_{i\ell}>0}
x_i^{b_{i\ell}}+ \mathop{\prod_{1\leq i\leq M}}_{b_{N+i,\ell}<0}
y_i^{-b_{N+i,\ell}} \mathop{\prod_{1\leq i\leq N}}_{b_{i\ell}<0}
x_i^{-b_{i\ell}}
\ee  %

Thus, the coefficients appearing in the exchange relation change with each
successive mutation, since they depend on the exchange matrix.

Let $\widetilde{\rho}=
\left( \begin{array}{cc} \rho & \mathbf{0} \\
\mathbf{0} & I_M \end{array} \right),$ where $\mathbf{0}$ denotes zeros and
$I_M$ denotes the $M\times M$ identity matrix. Thus $\widetilde{\rho}$
represents the permutation sending $(1,2,\ldots ,N)$ to $(N,1,2,\ldots ,N-1)$
and fixing $N+1,\ldots ,N+M$.

\bd
A quiver $Q$, with $N+M$ vertices as above, satisfying

\be \label{e:coefficientperiod1} %
\widetilde{\mu_1} B_Q=\widetilde{\rho} B_Q \widetilde{\rho}^{\, -1}
\ee  %
is said to be a \emph{period $1$ ice quiver}.
\ed

Since the effect of conjugation by $\widetilde{\rho}$ on the first $N$ elements
of each of the rows $N+1,\ldots ,N+M$ of the matrix $B_Q$ is to cyclically
shift them along one position to the right (with the entries in the opposite
positions in the extra columns cyclically moving down one position), it is easy
to see that, if we mutate such a quiver $Q$ successively at vertices
$1,2,\ldots ,N,1,2,\ldots $ etc. we obtain the terms (as in
Section~\ref{laurent}) of a recurrence:
$$
x_n x_{n+N} = F(x_{n+1},\ldots ,x_{n+N-1},y_1,y_2,\ldots ,y_M)
$$
where $F$ is a sum of two monomials in the $x_i$ with coefficients given by
fixed monomials in the $y_i$. By the Laurent Phenomenon~\cite[3.1]{02-3}, each
cluster variable can be written as a Laurent polynomial in $x_1,x_2,\ldots
,x_N$ with coefficients in $\mathbb{Z}[y_1,y_2,\ldots ,y_M]$. Thus, the
recurrence will be Laurent in this sense.

\bt \label{periodicicequivers}
Let $Q$ be an ice quiver on $N+M$ vertices, $1,2,\ldots ,N+M$, with vertices
$N+1,\ldots ,N+M$ frozen. Then $Q$ is a period $1$ ice quiver if and only
if the induced subquiver on vertices $1,2,\ldots ,N$ is a period $1$ quiver
and, if $B_Q$ is written as in Theorem~\ref{p1-prototheorem},
the following are satisfied.
\begin{description}
\item[(a)]
If $N=2r+1$ is odd, then for each $1\leq i\leq M$ such that row $N+i$ of
$B$ is non-zero there
is $t_i\in \{1,\ldots ,r\}$ such that $m_{t_i}=m_{N-t_i}=-1$ and all other $m_j$,
for $1\leq j\leq N-1$, are nonnegative, and a positive integer $l_i$ such that
$$b_{N+i,j}=\left\{ \begin{array}{cc}
l_i, & 1\leq j\leq t_i, \\
0, & t_i+1\leq j\leq 2r+1-t_i, \\
-l_i, & 2r+2-t_i\leq j\leq 2r+1=N, \\
0, & N+1 \leq j\leq N+M.
\end{array}\right.
$$
Alternatively, the $m_j$ are as above with the opposite signs and the entries
in the row are the negative of the above.
\item[(b)]
If $N=2r$ is even, then for each $1\leq i\leq M$ such that row $N+i$ of $B$
is non-zero there
is $t_i\in \{1,\ldots ,r-1\}$ such that $m_{t_i}=m_{N-t_i}=-1$ and the other
$m_j$, for $1\leq j\leq N-1$, are nonnegative, or $m_r=-2$ and
all other $m_j$ are nonnegative. Furthermore, there is a positive
integer $l_i$ such that
$$
b_{N+i,j}=\left\{ \begin{array}{cc}
l_i, & 1\leq j\leq t_i, \\
0, & t_i+1\leq j\leq 2r-t_i, \\
-l_i, & 2r+1-t_i\leq j\leq 2r=N, \\
0, & N+1\leq j\leq N+M.
\end{array} \right.
$$
Alternatively, the $m_j$ are as above with the opposite signs and the entries
in the row are the negative of the above.
\end{description}
\et

\noindent \textbf{Proof}:
To solve equation~(\ref{e:coefficientperiod1}), it is clear that the induced
subquiver of $Q$ on vertices $1,2,\ldots ,N$ must be a period $1$ quiver in our
usual sense. So we assume that the entries $b_{ij}$ for $1\leq i,j\leq N$ are
as in the general solution given in Theorem~\ref{p1-prototheorem}. For $1\leq
i\leq M$, the $N+i,j$ entry of $\widetilde{\rho} B_Q \widetilde{\rho}^{\, -1}$
is $b_{N+i,j-1}$ (where $j-1$ is read as $N$ if $j=1$). Thus we must solve the
equations $b_{N+i,N}=-b_{N+i,1}$ and
$b_{N+i,j-1}=b_{N+i,j}+\frac{1}{2}(b_{N+i,1}|b_{1,j}|+|b_{N+i,1}|b_{1,j})$ for
$i=1,2,\ldots ,M$ and $j=2,\ldots ,N$, noting that columns $N+1,\ldots ,N+M$
will give rise to the same equations and that $b_{N+i,j}=0$ for $j>N$.

We thus must solve the equations $b_{N+i,N}+b_{N+i,1}=0$ and
$b_{N+i,j-1}-b_{N+i,j}=\varepsilon(b_{N+i,1},m_{j-1})$ for $j=2,\ldots ,N$. Adding
all of these, we obtain the constraint that
\be  \label{c1-eq}  %
2b_{N+i,1}=\sum_{j=1}^{N-1} \varepsilon(b_{N+i,1},m_j).
\ee  %
The solutions are given by the values of $b_{N+i,1}$ satisfying this constraint
(using the other equations to write down the values of the other $b_{N+i,j}$).

If $N=2r+1$ is odd, the constraint can be rewritten
$$
2b_{N+i,1}=b_{N+i,1}(|m_1|+\cdots +|m_r|)-|b_{N+i,1}|(m_1+\cdots +m_r),
$$
using the fact that $m_j=m_{N-j}$ for all $j$. If $b_{N+i,1}=0$ then
$b_{N+i,j}=0$ for all $j$, so for a non-zero $(N+i)$th
row we must have $b_{N+i,1}\not=0$. For solutions with $b_{N+i,1}>0$, we must have
$$
|m_1|-m_1+|m_2|-m_2+\cdots + |m_r|-m_r=2.
$$
Since $|x|-x=0$ for $x\geq 0$ and equals $-2x$ for $x\leq 0$, the only
solutions arise when $m_{t_i}=-1$ for some $t_i$ (so $m_{N-t_i}=-1$ also)
and all other $m_j$ are nonnegative.
They are of the form
$$
b_{N+i,j}=\left\{ \begin{array}{cc}
b_{N+i,1}, & 1\leq j\leq t_i \\
0, & t_i+1\leq j\leq 2r+1-t_i \\
-b_{N+i,1}, & 2r+2-t_i\leq j\leq 2r+1,
\end{array}\right.
$$
as required.
The solutions with negative $b_{N+i,1}$ are the negative of these
(provided $m_{t_i}=m_{N-t_i}=1$ and all other $m_j$ are nonpositive).

If $N=2r$ is even, the constraint can be rewritten
$$
2b_{N+i,1}=b_{N+i,1}(|m_1|+\cdots +|m_{r-1}|+|m_r|/2)-|b_{N+i,1}|(m_1+\cdots +m_{r-1}+m_r/2),
$$
again using the fact that $m_j=m_{N-j}$ for all $j$. As in the odd case,
we must have $b_{N+i,1}\not=0$ for a non-zero $(N+i)$th row.
If $b_{N+i,1}>0$, we must have
$$
|m_1|-m_1+|m_2|-m_2+\cdots + |m_{r-1}|-m_{r-1}+|m_r|/2-m_r/2=2.
$$
Arguing as in the odd case, we see that solutions arise when $m_{t_i}=-1$ for
some $t_i$ with $1\leq t_i\leq r-1$ (and so $m_{N-t_i}=-1$) and all other
$m_j$ are nonnegative, or when $m_r=-2$ and all other $m_j$ are nonnegative.
They are of the form
$$
b_{N+i,j}=\left\{ \begin{array}{cc}
b_{N+i,1}, & 1\leq j\leq t_i \\
0, & t_i+1\leq j\leq 2r-t_i \\
-b_{N+i,1}, & 2r+1-t_i\leq j\leq 2r,
\end{array} \right.
$$
where $t_i=r$ for the last case, as required.
The solutions for negative $b_{N+i,1}$ are the
negative of these (with the negative of the constraints on the $m_j$).
$\Box$

\bc \label{othermonomial}
Consider the Laurent recurrence~\eqref{period1recurrence} corresponding to
a period $1$ quiver:
$$x_nx_{n+N}=\prod_{\substack{i=1\\ m_i>0}}^{N-1}x_{n+i}^{m_i}+
\prod_{\substack{i=1\\ m_i<0}}^{N-1}x_{n+i}^{-m_i}$$
The same recurrence, with parameters introduced on the right hand
side as coefficients of the monomials, arises from a period $1$
ice quiver as above if and only if a parameter on a monomial is
only allowed when the \emph{other} monomial is of the form
$x_{n+i}x_{n+N-i}$ for some $i$ with $1\leq i\leq N/2$.
(If both monomials are of this form then a parameter is allowed on both.)
If this condition is satisfied, the recurrence with parameters is again Laurent.
\ec

This corollary has the following interesting consequence:

\bp \emph{\textbf{(Gale-Robinson recurrence)}}
The only binomial recurrences corresponding to period $1$ quivers
that, when parameters on both monomials are allowed, correspond to period
$1$ ice quivers, are the two-term Gale-Robinson recurrences.
\ep

\noindent \textbf{Proof}:
The Gale-Robinson recurrences are exactly those for which both monomials
are of the required form in Corollary~\ref{othermonomial}. See
Example~\ref{galerob} in Section~\ref{p1-sequence}.
$\Box$

Note that it follows that these recurrences, with a parameter multiplying each
of the monomials, are Laurent. This was shown in~\cite[1.7]{02-2}.

\begin{figure}[ht]
\centering
\includegraphics[width=4.5cm]{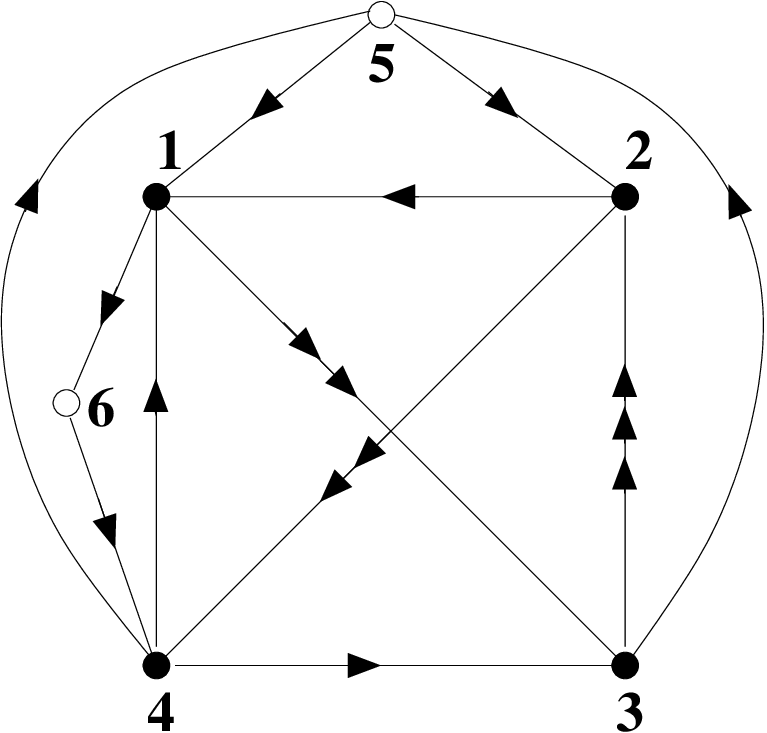}
\caption{The ice quiver for Somos $4$ with
parameters.}\label{fig:somos4parameterquiver}
\end{figure}

\bex[Somos $4$ Recurrence with Parameters]  \label{ex:somos4parameters} {\em
The Somos $4$ recurrence is a special case of the two-term Gale-Robinson
recurrence. We can add extra coefficient rows $(1,1,-1,-1,0,0)$ and
$(-1,0,0,1,0,0)$ to the corresponding matrix, giving the six-vertex quiver
shown in Figure~\ref{fig:somos4parameterquiver}, with empty circles denoting
frozen vertices. We recover the Laurent property of the corresponding
recurrence:

$$
x_{n+4}x_n=y_1x_{n+1}x_{n+3}+y_2x_{n+2}^2.
$$
}\eex

\bex[A Recurrence Considered by Dana Scott] {\em We have seen (see
Example~\ref{4node} of Section \ref{p1-sequence}) that the $4$-node case
$m_1=1$, $m_2=-1$ corresponds to a recurrence considered by Dana Scott.
By Theorem~\ref{periodicicequivers}(b) the only possible non-zero extra
row in the matrix is $(-1,0,0,1,0,0)$ (or a positive multiple) giving the
recurrence:
$$
x_{n+4}x_n=x_{n+1}x_{n+3}+yx_{n+2}.
$$
It follows that this recurrence is Laurent. The corollary, that this recurrence
gives integers for all integer $y$ (if $x_1=x_2=x_3=x_4=1$), was noted
in~\cite{91-17}. It is also noted in~\cite{91-17} that the recurrence
$$
x_{n+4}x_n=2x_{n+1}x_{n+3}+x_{n+2}
$$
(with the same initial conditions) does not give integers. It follows that the
recurrence
$$
x_{n+4}x_n=y_1x_{n+1}x_{n+3}+x_{n+2}
$$
does not have the Laurent property.
}\eex   %

\br[A Conjecture] %
Although our results do not determine non-Laurentness it is interesting to note
that we obtain exactly those which are Laurent as solutions in the above
examples. It seems reasonable to conjecture that the parameter versions of
recurrences arising from period $1$ quivers are Laurent if and only if they
arise from a cluster algebra with frozen variables in the above sense.
\er   %

\section{Supersymmetric Quiver Gauge Theories}\label{susy}

\begin{figure}[ht]
\centering \subfigure[Hirzebruch $0$]{
\includegraphics[width=2.5cm]{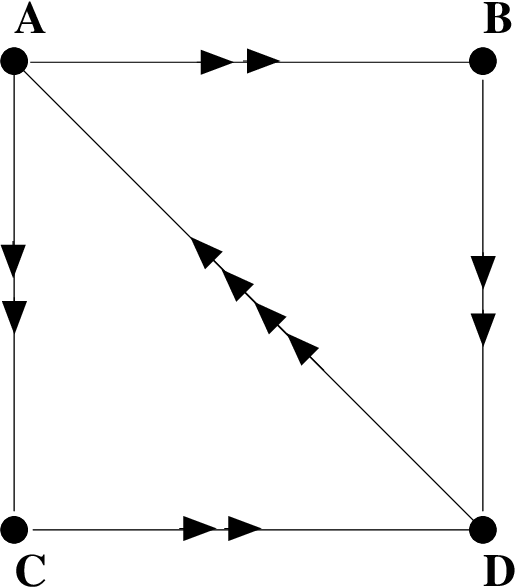}\label{fig:hirzebruch0}
} \subfigure[del Pezzo $0$]{
\includegraphics[width=2.5cm]{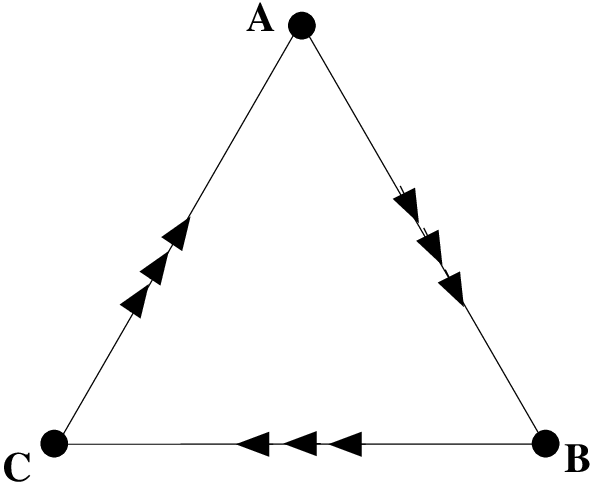}\label{fig:delPezzo0}
} \subfigure[del Pezzo $1$]{
\includegraphics[width=2.5cm]{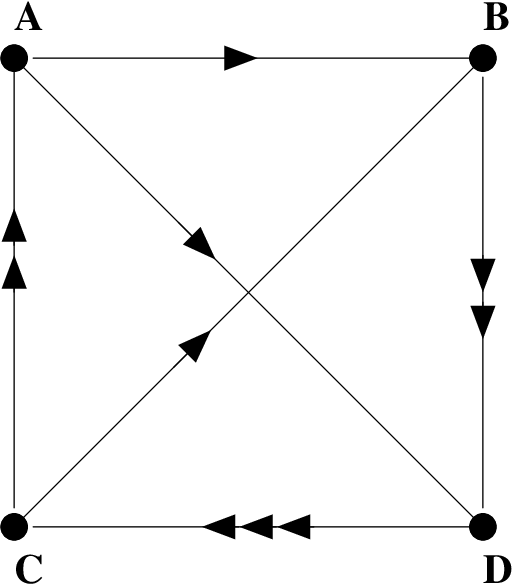}\label{fig:delPezzo1}
}
\subfigure[del Pezzo $2$]{
\includegraphics[width=2.5cm]{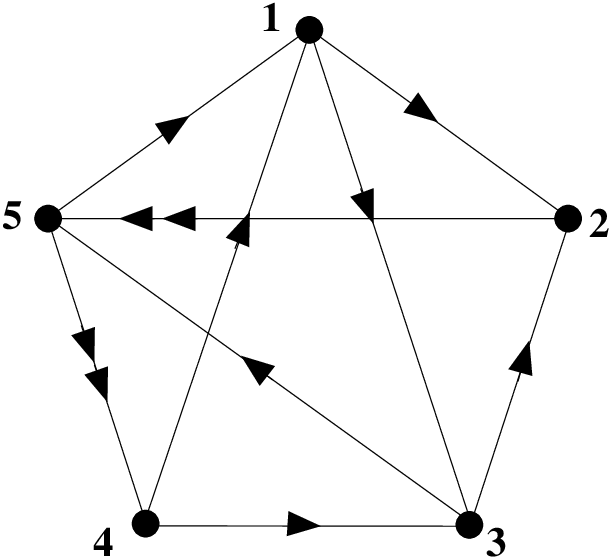}\label{fig:MRdelPezzo2}
} \subfigure[del Pezzo $3$]{
\includegraphics[width=3cm]{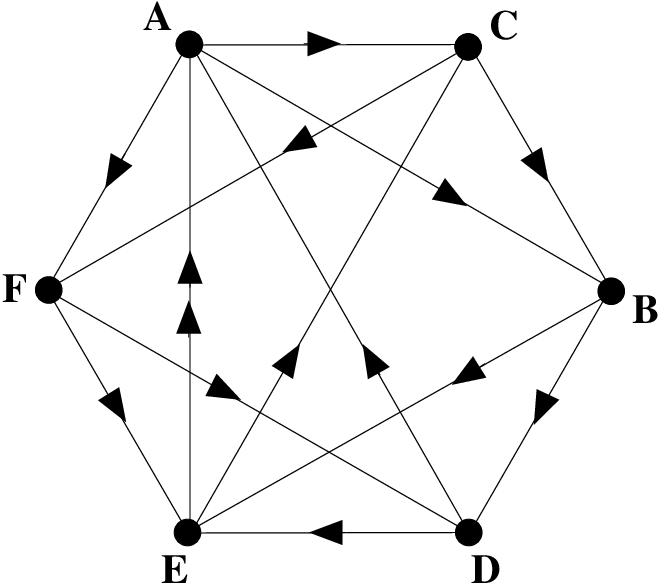}\label{fig:delPezzo3}
} \label{fhhquivers} \caption{Quivers of quiver gauge theories associated to a
family of surfaces. Quiver~\subref{fig:MRdelPezzo2} is from~\cite{04-3} while
the others are from~\cite{01-10}.}
\end{figure}

In this section we point to the \emph{$D-$brane} literature in which our
quivers arise in the context of \emph{quiver gauge theories}.
The quivers arising in supersymmetric quiver gauge theories often have
periodicity properties. Indeed, in~\cite[\S3]{04-3} the
authors
consider an $\mathcal{N}=1$ supersymmetric quiver gauge theory associated
to the complex cone over the second del Pezzo surface $dP_2$. The quiver
$Q_{dP_2}$ of the gauge theory they consider is given in
Figure~\ref{fig:MRdelPezzo2}. The authors compute the Seiberg dual
of the quiver gauge theory at each of the nodes of the quiver. The Seiberg
dual theory has a new quiver, obtained using a combinatorial rule from
the original quiver using the choice of vertex (see~\cite{01-9}).
It can be checked that the combinatorial rule for Seiberg-dualising a
quiver coincides with the rule for Fomin-Zelevinsky quiver mutation
(Definition~\ref{d:mutate}); see~\cite{07-4} for a discussion of
the relationship between Seiberg duality and quiver mutation.
In~\cite[\S 3]{04-3} the authors compute the Seiberg dual of
$Q_{dP_2}$ at each node, in particular showing that the Seiberg dual of
$Q_{dP_2}$ at node $1$ is an isomorphic quiver. They indicate that such
behaviour is to be expected from a physical perspective.

This quiver fits into the scheme discussed in this article: it is a period $1$
quiver. In fact it coincides with the quiver corresponding to the matrix
$B_5^{(1)}-B_5^{(2)}+ B_3^{(1)}$ (see Example~\ref{n=5p1} with $m_1=1,m_2=-1$,
and also Figure~\ref{fig:somos5quiver}), with the relabelling
$(1,2,3,4,5)\mapsto (3,4,5,2,1)$ (we give all relabellings starting from our
labels). We note that this quiver also appears in \cite{07-3} (with a
relabelling $(1,2,3,4,5)\mapsto (2',3',1,2,3)$) in the context of a $dP_2$
brane tiling and that the corresponding sequence is the Somos $5$ sequence.

The quivers of quiver gauge theories associated to the complex cones over the
Hirzebruch zero and del Pezzo $0$--$3$ surfaces are computed in~\cite[\S
4]{01-10}. We list them for convenience in Figure~\ref{fhhquivers} for the
Hirzebruch $0$ and del Pezzo $0,1$ and $3$ surfaces. Note that the del Pezzo
$2$ case was discussed above: we chose the quiver given in~\cite[\S3]{04-3} for
this case because it fits better into our setup.  We remark that the quiver for
$dP_1$ coincides with the Somos $4$ quiver (with matrix
$B_4^{(1)}-2B_4^{(2)}+2B_2^{(1)}$), with the relabelling $(1,2,3,4)\mapsto
(B,C,D,A)$. Thus it is period $1$. See Example~\ref{n=4p1} and
Figure~\ref{subfig:somos4quiver}.

The quiver for the Hirzebruch $0$ surface is period $2$: with the relabelling
$(1,2,3,4)\mapsto (C,A,D,B)$. It corresponds to the matrix $B(1)$ given in
equation~(\ref{p2n4}) in Section~\ref{p2-4node} with $m_1=2$, $m_2=-2$ and
$m_3=0$. The quiver for $dP_3$ is period two. In fact it is one of the period
two quivers described in Section~\ref{p2-reg}, with $m_1=m_3=1$, $m_2=-1$ and
$m_{\overline{1}}=0$.

Finally, we note that, by construction, the in-degree of a vertex $i$ always
coincides with the out-degree of $i$ for any quiver arising from a brane tiling
in the sense of \cite{06-3}. It is interesting to note that the only quivers of
cluster mutation period $1$ satisfying this assumption with five or fewer
vertices are the Somos 4 and Somos 5 quivers, i.e. quivers associated to $dP_1$
and $dP_2$ (see Figures \ref{fig:delPezzo1} and \ref{fig:MRdelPezzo2}).

\section{Conclusions}\label{conclude}

In this paper we have raised the problem of classifying all quivers with
mutation periodicity. For period 1 we have given a {\em complete} solution. For
period 2 we have given a solution which exists {\em for all $N$} (the number of
nodes).  In addition to these, we have seen that, for $N=5$ there are some
``exceptional solutions'', such as (\ref{n5-1pos4neg2pos}) and
(\ref{n5-1pos4neg2neg}).  We conjecture that such ``exceptional solutions''
will exist for all {\em odd $N$}, but not for the even case.  We also
conjecture that there are more general infinite families than the one presented
in Section \ref{p2-reg}, where we took a particularly simple initial condition
for the iteration (\ref{eqn:period2equation}).  We could, for instance,
consider the case with $m_{N-r}=m_{\overline{r}}$ for several values of $r$. We
could also consider a family of {\em period $r$} quivers, satisfying
$\mu_1(B)=\rho B(\sigma(\mathbf{m}))\rho^{-1}$, with
$\sigma^r(\mathbf{m})=\mathbf{m}$.

The other main theme of our paper was the construction and analysis of recurrences
with the Laurent property.  We have shown that the recurrences associated with our
period 1 primitives can be linearised.  This construction can be extended to
higher period cases, but we currently only have examples.  General periodic
quivers give rise to truly nonlinear maps, the simplest of which is Somos $4$,
which is known \cite{07-2} to be integrable and, in fact, related to the QRT
\cite{88-5} map of the plane.  Somos $5$ is similarly known to be integrable,
as are the subcases of Somos $6$ discussed in Example \ref{6node1}. On the
other hand,
$$
x_nx_{n+6} = x_{n+1}^2x_{n+5}^2+x_{n+2}^2x_{n+3}^4x_{n+4}^2,
$$
(corresponding to the choice $m_1=-m_2=2, m_3=-4$ in Example \ref{6node1}) is
known to be {\em not} integrable (see Equation (4.3) of \cite{07-2}). Even
though this recurrence has the Laurent property and satisfies ``singularity
confinement'' \cite{91-6} (a type of Painlev\'e property for discrete
equations), it fails the more stringent ``algebraic entropy'' \cite{99-12} test
for integrability. This simple test (or the related ``diophantine
integrability'' \cite{05-4} test) can very quickly show that a map is {\em not}
integrable. If they {\em indicate} integrability, then it is sensible to search
for {\em invariant functions} in order to {\em prove} integrability. Early
indications (preliminary calculations by C-M Viallet) are that for integers
$m_i$ which are ``small in absolute value'' we have integrable cases. We thus
expect that small sub-families of our general periodic quivers will give rise
to integrable maps.
The isolation and classification of integrable cases of the recurrences discussed in this
article is an important open question, which will be discussed in~\cite{FordyHone}.

We have seen that many of our examples occur in the context of supersymmetric
quiver gauge theories.  A deep understanding of the connection with brane
tilings and related topics would be very interesting.

After this paper first appeared on the arXiv, Jan Stienstra pointed out to us
that the quivers in Figure~\ref{fhhquivers} also appear in~\cite{08-4} in the
context of Gelfand-Kapranov-Zelevinsky hypergeometric systems in two variables,
suggesting a possible connection between cluster mutation and such systems.

Since we wrote the first version, the article~\cite{nakanishi} has appeared,
proposing a general study of periodicity in cluster algebras (in a wide
sense), motivated by many interesting examples of periodicity for cluster algebras
and $T$ and $Y$ systems (see references therein).
In particular, periodicity of (seeds in) cluster algebras plays a key role
in the context of the periodicity conjecture for $Y$ systems, which was proved
in full generality in~\cite{keller}; see also references therein.

The articles~\cite{assemreutenauersmith,kellerscherotzke} have also appeared,
proving the linearisation of frieze sequences (or frises) associated to acyclic
quivers (quivers without an oriented cycle); see these two papers and
references therein for more details of these sequences. If $Q$ is acyclic, its
vertices $\{1,2,\ldots ,N\}$ can be numbered so that $i$ is a sink on the
induced subquiver on vertices $i,i+1,\ldots ,N$ for each $i$; then $1,2,\ldots ,N$ is
an admissible sequence of sinks in the sense of Remark~\ref{r:admissible} and so
$Q$ has period $N$ in our sense. The corresponding sequence of cluster variables can
be regarded as an $N$-dimensional recurrence (as in the period $2$
case,~\ref{ss:period2case}). Generalising results of~\cite{assemreutenauersmith},
in~\cite{kellerscherotzke} it is shown that all components of this recurrence
are linearisable if $Q$ is Dynkin or affine (and conversely; in fact the result is more general,
including valued quivers). This linearisability can be regarded as a generalisation of
Theorem~\ref{linear-thm}, since our primitives are acyclic and the numbering of their vertices
satisfies the above sink requirement, using Lemma~\ref{taupreserves}. In fact this last statement
is true for primitives of any period, so it follows from~\cite{kellerscherotzke}
that the recurrences corresponding to primitives of any period are linearisable, noting that the period of a
primitive always divides $N$.

\subsection*{Acknowledgements}

We would like to thank Andy Hone for some helpful comments on an earlier version of
this manuscript, Jeanne Scott for some helpful discussions, and Claude Viallet for computing the
algebraic entropy of some of our maps.
We'd also like to thank the referees for their helpful and interesting
comments on the submitted version of this article. This work was
supported by the Engineering and Physical Sciences Research Council [grant
numbers EP/C01040X/2 and EP/G007497/1].

\section{Appendix A1: Proofs of Results in Section
\ref{lin:diff:eq}}\label{append}

\noindent \textbf{Proof of Lemma~\ref{lem:andescription}}: We first prove (a)
and (b). The result is clearly true for $n=0,1$. Assume that it holds for
smaller $n$ and firstly assume that $n=2r-1$ is odd. Then
\begin{align*}
a_{2r-1} &= -a_{2r-3}-c_{2r-2}a_{2r-2} \\
&= (-1)^{r-1} \sum_{k=0}^{r-2} (-1)^k t_{2k,\text{odd}}^{(2r-4)}
+c_{2r-2} (-1)^{r} \sum_{k=0}^{r-2} (-1)^k t_{2k+1,\text{odd}}^{(2r-3)} \\
&= (-1)^{r-1}t_{0,\text{odd}}^{(2r-4)}+ (-1)^{r-1} \sum_{k=0}^{r-3} (-1)^{k+1}
t_{2k+2,\text{odd}}^{(2r-4)}
+c_{2r-2} (-1)^{r} \sum_{k=0}^{r-2} (-1)^k t_{2k+1,\text{odd}}^{(2r-3)} \\
&= (-1)^{r-1}t_{0,\text{odd}}^{(2r-4)}+
(-1)^{2r-2}c_{2r-2}t_{2r-3,\text{odd}}^{(2r-3)}+ (-1)^r \sum_{k=0}^{r-3} (-1)^k
(t_{2k+2,\text{odd}}^{(2r-4)}+
c_{2r-2}t_{2k+1,\text{odd}}^{(2r-3)}) \\
&= (-1)^{r-1}t_{0,\text{odd}}^{(2r-4)}+
(-1)^{2r-2}t_{2r-2,\text{odd}}^{(2r-2)}+
(-1)^r \sum_{k=0}^{r-3} (-1)^k t_{2k+2,\text{odd}}^{(2r-2)} \\
&= (-1)^{r-1} \sum_{k=0}^{r-1} (-1)^k t_{2k,\text{odd}}^{(2r-2)},
\end{align*}
and the result holds for $n$. A similar argument shows that the result holds
for $n$ when $n$ is even. Then (a) and (b) follow by induction. To prove (c),
we note that $t_{2k+1,\text{odd}}^{(2r-1)}$ is invariant under the
transformation $c_l\mapsto c_{2r-l}$ and that $t_{2k,\text{odd}}^{(2r-2)}$ is
invariant under $c_l\mapsto c_{2r-1-l}$. We then have
$$b_{N-1}=a_{N-1}|_{c_l\mapsto c_{2N-2-l}}=a_{N-1}|_{c_l\mapsto c_{N-1-l}}=a_{N-1}$$
using Lemma~\ref{lem:abflip} and the fact that $c_{n+N-1}=c_n$.$\Box$

\medskip\noindent \textbf{Proof of Proposition~\ref{prop:linearisation}}: We
first assume that $N=2r-1$ is odd, so $(-1)^{N-2}=-1$. Then
{\allowdisplaybreaks
\begin{align*}
S_{N,1} &= -(a_{2r-3}+a_{2r-3}|_{c_l\mapsto 2r-2-l}+c_{2r-2}a_{2r-2}) \\
&= (-1)^{r-1} \sum_{k=0}^{r-2} (-1)^k t_{2k,\text{odd}}^{(2r-4)} +(-1)^{r-1}
\sum_{k=0}^{r-2} (-1)^k t_{2k,\text{even}}^{(2r-3)}
+c_{2r-2}(-1)^{r}\sum_{k=0}^{r-2}(-1)^k t_{2k+1,\text{odd}}^{(2r-3)} \\
&= (-1)^{r-1}t_{0,\text{odd}}^{(2r-4)}+ (-1)^{r-1} \sum_{k=0}^{r-3} (-1)^{k+1}
t_{2k+2,\text{odd}}^{(2r-4)}
+(-1)^{r-1} \sum_{k=0}^{r-2} (-1)^k t_{2k,\text{even}}^{(2r-3)} \\
& \quad\quad
+c_{2r-2}(-1)^{r}\sum_{k=0}^{r-2}(-1)^k t_{2k+1,\text{odd}}^{(2r-3)} \\
&= (-1)^{r-1}t_{0,\text{odd}}^{(2r-4)}+ c_{2r-2}(-1)^{2r-2}
t_{2r-3,\text{odd}}^{(2r-3)}+ (-1)^{r} \sum_{k=0}^{r-3} (-1)^{k}
(t_{2k+2,\text{odd}}^{(2r-4)}+
c_{2r-2}t_{2k+1,\text{odd}}^{(2r-3)}) \\
& \quad\quad +(-1)^{r-1} \sum_{k=0}^{r-2} (-1)^k t_{2k,\text{even}}^{(2r-3)}
\\
&= (-1)^{r-1}t_{0,\text{odd}}^{(2r-2)} +(-1)^{2r-2}
t_{2r-2,\text{odd}}^{(2r-2)}+ (-1)^{r} \sum_{k=0}^{r-3} (-1)^{k}
t_{2k+2,\text{odd}}^{(2r-2)}
+(-1)^{r-1} \sum_{k=0}^{r-2} (-1)^k t_{2k,\text{even}}^{(2r-3)} \\
&= (-1)^{r-1} \sum_{k=0}^{r-1} (-1)^{k} t_{2k,\text{odd}}^{(2r-2)}
+(-1)^{r-1} \sum_{k=0}^{r-2} (-1)^k t_{2k,\text{even}}^{(2r-2)}\\
&= (-1)^{r-1} \sum_{k=0}^{r-1} (-1)^{k} t_{2k,\text{alt}}^{(2r-2)},
\end{align*}}
as required. A similar argument shows that the result holds when $N$ is even.

As we have already observed, for the sequence $(x_i)$ we are interested in, the
$L_i$ vanish. It follows that $L$ vanishes and we are done.$\Box$

\medskip\noindent \textbf{Proof of Corollary~\ref{cor:generallinearisation}}:
The proof of Proposition~\ref{prop:linearisation} also shows that
$x_2+x_{2+2N-1}=S_{N,1}|_{c_l\mapsto c_{l+1}}x_{2+N-1}$. It follows from the
description of $S_{N,1}$ in Proposition~\ref{prop:linearisation} that
$S_{N,1}=S_{N,1}|_{c_l\mapsto c_{l+1}}$ (using the fact that $c_{n+N-1}=c_n$)
so we are done for $n=2$. Repeated application of this argument gives the
result for arbitrary $n$.$\Box$

\end{document}